\documentclass[12pt,reqno]{paper}
\usepackage{fancyhdr}
\usepackage[a4paper, hmargin=2.5cm, vmargin={3.2cm, 3.2cm}]{geometry}

\usepackage[pdftex]{graphicx}
\usepackage{titlesec}
\titleformat{\subsection}[runin]
{\normalfont\bfseries}{\thesubsection}{1em}{}
\titleformat{\subsubsection}[runin]
{\normalfont\bfseries}{\thesubsubsection}{1em}{}

\usepackage[titletoc,title]{appendix}

\usepackage{enumerate}
\usepackage{amsmath}%
\usepackage{amsfonts}%
\usepackage{amsthm}
\usepackage{amssymb}%
\usepackage{graphicx}
\usepackage{theoremref}
\usepackage{dsfont}

\DeclareMathOperator{\re}{Re}
\DeclareMathOperator{\im}{Im}

\theoremstyle{plain}

\numberwithin{equation}{subsection}

\newcommand{\e}{\varepsilon}
\newcommand{\ph}{\varphi}

\newtheorem{theorem}{Theorem}

\newtheorem{lemma}{Lemma}
\numberwithin{lemma}{subsection}
\newtheorem*{Admissible}{Definition}

\title{On functions $K$ and $E$ generated by a sequence of moments}
\author{Avner Kiro
\thanks{Both authors are partially supported by ISF Grants~166/11 and~382/15 and BSF Grant~2012037.}
\\
School of Mathematical Sciences, Tel Aviv University\\
Tel Aviv 69978, Israel
\\
E-mail address: {\tt avnerefrak@mail.tau.ac.il}
\and
Mikhail Sodin
\\
School of Mathematical Sciences, Tel Aviv University\\
Tel Aviv 69978, Israel
\\
E-mail address: {\tt sodin@post.tau.ac.il}
}

\begin{document}
\maketitle

\begin{abstract}
We study the asymptotic behaviour of the entire function
\[
E(z) = \sum_{n\ge 0} \frac{z^n}{\gamma (n+1)}
\]
and the analytic function
\[
K(z) = \frac1{2\pi {\rm i}}\, \int_{c-{\rm i}\infty}^{c+{\rm i}\infty} z^{-s}\gamma (s)\, {\rm d}s\,,
\]
which naturally appear in various classical problems of analysis.
\end{abstract}

\section{Introduction and main results}

\subsection{The functions $K$ and $E$.}

In this work we study the asymptotic behavior of two analytic functions $K$ and $E$ generated by a sequence of moments $(\gamma(n+1))_{n\geq0}$, where $\gamma(s)$ is an analytic function in the angle $\{s: |\arg(s+c)|<\alpha_0\}$ with $\frac{\pi}{2}<\alpha_0\leq\pi$ and $c=c_\gamma>0$.
The function $\gamma$ satisfies certain regularity properties, which we will list shortly. Here, we will only mention that $(\gamma(n))$ is a fastly growing sequence of positive numbers (so that $\displaystyle \lim_{n\to\infty}\gamma(n)^{1/n}=\infty$), and that, for some $\alpha\in(\frac{\pi}{2},\alpha_0)$,
\[ \lim_{\rho\to\infty} \frac{\log|\gamma(\rho e^{\pm i\alpha})|}{\rho}=-\infty\,. \]
This allows us to define the functions
\begin{equation}
 K(z)=\frac{1}{2\pi i}\int_{\mathcal{L}_\alpha} z^{-s}\gamma(s)\, {\rm d}s\,, \label{Kdef}
\end{equation}
where $\mathcal{L}_\alpha=\{z : |\arg(z)|=\alpha\}$ is a union of two rays traversed
is such a way that $\im(z)$ increases along  $\mathcal{L}_\alpha$ (see Figure 1), and
\begin{equation}
E(z)=\sum_{n\geq 0} \frac{z^n}{\gamma(n+1)}\,. \label{Edef}
\end{equation}
\begin{figure}[h]
	\centering
	\includegraphics[scale=0.35]{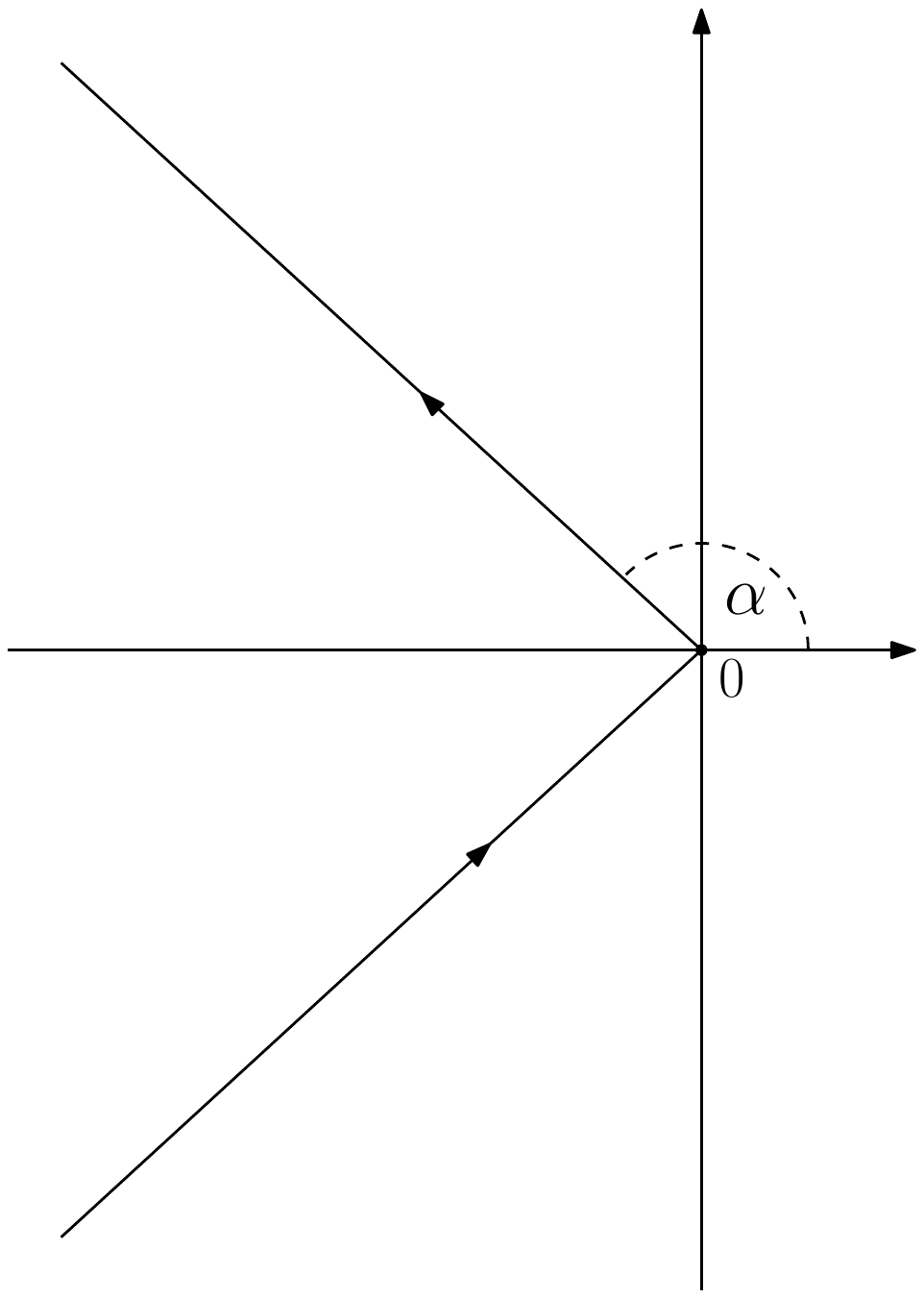}
	\caption{Contour  $\mathcal{L}_\alpha$}
\end{figure}
The function $K$ is analytic on the Riemann surface of $\log z$ (that is, the function $K(e^w)$ is entire), while the function $E$ is an entire one.

The assumptions on the function $\gamma$, which we will impose shortly, will allow us, moving the integration contour, represent the function $K(t)$ for $t\geq 0$ as the inverse Mellin transform of $\gamma$:
\begin{equation}\label{eq:K_Mellin}
K(t)=\frac{1}{2\pi {\rm i}}\int_{c-{\rm i}\infty}^{c+{\rm i}\infty}t^{-s}\gamma(s)\,
{\rm d}s,\quad c>0\,,
\end{equation}
where the integral does not depend on the choice of $c>0$. Then, by the inversion formula for the Mellin transform, $K$ solves the moment problem
\[
\int_0^\infty t^n K(t)\, {\rm d}t = \gamma(n+1),\quad n\in \mathbb{Z}_+\,.
\]

The functions $K$ and $E$ naturally appear in various classical problems of analysis, for instance, in the Borel-type moment summation of divergent series \cite{Hardy} and in studying in convergence of certain interpolation problems for entire functions \cite{Evgrafov1, Evgrafov2, Evgrafov3}. It also worth mentioning 
that Beurling~\cite{Beurling2, Beurling}
singled out a class of functions $\gamma$ for which the function $K$ is positive on the
positive half-line (see Section~1.7 and Appendix~B). Then, our Theorem~1 gives explicit
asymptotics of solutions to a large class of determinate Stieltjes moment problem.

Our interest originated in Beurling's approach to the problems of description of the Taylor coefficients and of summation of the divergent Taylor series in various classes of smooth functions \cite{Beurling2, Beurling, Kiro}.

It could be that the results presented here are known to experts. On the other hand, we were unable to locate them in the literature and we believe that they are of certain interest.

\medskip
At last, we note that juxtaposing the right-hand sides of \eqref{Kdef} and \eqref{Edef}, we may expect some match between the growth of $E$ and the decay of $K$ on the positive half-line and nearby. In the prototype case when $\gamma$ is Euler's gamma--function, $\gamma(n+1)=n!$, $E(z)=e^z$, $K(z) =e^{-z}$, and this match is prefect.

\subsection{Admissible functions.}
From now on, we will assume that the function $\gamma$ is analytic and non-vanishing in the angle
$\{s: |\arg(s+c)|<\alpha_0\}$ with $\frac{\pi}{2}<\alpha_0\leq\pi$ and $c=c_\gamma>0$,
and is positive on $(-c, +\infty)$. We put
\[
L(s)=\gamma(s)^{1/s} \quad \text{and} \quad \varepsilon(s)=s\,\frac{L'(s)}{L(s)}\,.
\]
\begin{Admissible}
	We call the function $\gamma$ {\em admissible}, if the function $\varepsilon$ is positive
and bounded on $\mathbb{R}_+$, and satisfies the following conditions:
	\begin{itemize}
		\item[{\rm (A)}]
$ \displaystyle \int^\infty \frac{\varepsilon(\rho)}{\rho}\, {\rm d}\rho = \infty$,
        \item [{\rm (B)}]$\rho|\varepsilon'(\rho)|=o\left(\varepsilon(\rho)\right)$ as $\rho\to\infty$,
		\item [{\rm (C)}] for $s=\rho e^{i\theta} $, $\rho\to\infty$,
		one has $\varepsilon (s)=\left(1+o(1)\right)\varepsilon(\rho)$,
		uniformly in any angle $|\theta|\leq \alpha_0-\delta$.
	\end{itemize}
	\end{Admissible}
\noindent Condition (A) means that the function $L$ is unbounded. Condition (B) says that the
function $\varepsilon$ is slowly varying.
Everywhere below, {\em we always assume that the function $\gamma$ is admissible}.

\medskip
	It is not difficult to see that conditions (B) and (C) yield that, for $s=\rho e^{i\theta}$, $\rho\to\infty$,
\begin{itemize}
\item[(D)] $\displaystyle \log L\left(\rho e^{i\theta}\right)=\int_{0}^{\rho}\frac{\varepsilon\left(u\right)}{u}\, {\rm d}u
    +i\theta\varepsilon(\rho)+o\left(\varepsilon\left(\rho\right)\right)$,
\item[(E)] $s\varepsilon'(s)=o(\varepsilon(\rho))$,
    \end{itemize}
also uniformly in any angle $|\theta|\le\alpha_0-\delta$.

\medskip
Indeed, condition (D) follows from (C) by integration, while (E) follows from (C) and (B)
due to the analyticity of $\varepsilon$.

\medskip
Below, in Sections~1.5 and 1.6, we will give several examples and constructions of admissible functions $\gamma$.

\subsection{The saddle-point equation.} It is clear, at least intuitively, that the asymptotic behavior of the functions $K(z)$ and $E(z)$ for large $z$ should be determined by the saddle-point of the function $s\mapsto \log \gamma(s)-s\log z=s\log L(s)-s\log z$, that is by the equation
\begin{equation}
\log L(s)+s\frac{L'(s)}{L(s)}=\log z,\label{saddlePoint}
\end{equation}
which we will call \textit{the saddle-point equation}.

For $0<\alpha<\alpha_0$ and $\rho_0>0$, put
\[ S(\alpha, \rho_0)=\{s: \; |\arg(s)|<\alpha,\; |s|>\rho_0\} \]
and let
\[
\Phi (s) = \log L(s) + s\,\frac{L'(s)}{L(s)}\,.
\]
Note that for $s\in S(\alpha, \rho_0)$ and $|s|=\rho$ sufficiently large, by (B),(C) and (E),
we have
\[ \re\Phi'(s) = \left(1+o(1)\right)\frac{\varepsilon(\rho)}{\rho}>0\,. \]
Thus, for $s_1, s_2\in S(\alpha, \rho_0)$, $s_1\ne s_2$,
\[
\re\,\frac{\Phi(s_2)-\Phi(s_1)}{s_2-s_1}
= \int_0^1 \re \Phi'(s_1+t(s_2-s_1))\, {\rm d}t >0\,,
\]
provided that $\rho_0$ is large enough.
Therefore, for $\rho_0$ sufficiently  large, the function $\Phi$, that is, the LHS of the saddle-point equation, is a univalent function in $S(\alpha, \rho_0)$.  From here on, we assume that this is the case. Then, we put
\[ \Omega(\alpha, \rho_0)=\left\{\log z= \Phi(s)\colon s\in S(\alpha, \rho_0)\right\}. \]
This is a domain on the Riemann surface of $\log z$. If the index $\rho_0$ is not essential,
we will skip it, to simplify notation.

Note that if $\rho=|s|$ is sufficiently large, by (C) and (D), we have
\[ \im \Phi (s) =\left(\theta+o(1)\right)\varepsilon(\rho),\quad s=\rho e^{i\theta}. \]
Thus, choosing $\rho_0$ sufficiently large, we can treat  $\Omega(\alpha, \rho_0)$ as a subdomain 
of the slit plane $\mathbb{C}\setminus\mathbb{R}_-$, provided that
\[ \limsup_{\rho\to\infty} \varepsilon(\rho)<\frac{\pi}{\alpha}, \]
in particular, whenever $\varepsilon(\rho)=o(1)$, as $\rho\to\infty$.

\medskip
In what follows,
by $s_z=\rho_z e^{{\rm i}\theta_z}$  we always denote the unique solution of the saddle-point equation \eqref{saddlePoint}.

\subsection{Asymptotics of the functions $K$ and $E$.}

We are now able to present our results.

\begin{theorem}\label{TheoremK}
Suppose that the function $\gamma$ is admissible. Then, for any $\delta>0$,  we have
\[K(z)=\left(1+o(1)\right)\sqrt{\frac{L(s)}{2\pi L'(s)}}\exp\left(-s^2\frac{L'(s)}{L(s)}\right), \quad z\to\infty\]
uniformly in $\Omega(\alpha_0-\delta)$. Here $s=s_z$ and the   branch of the  square root is  positive on the positive half-line.
\end{theorem}

\begin{theorem}\label{TheoremE}
Suppose that the function $\gamma$ is admissible and that
\begin{equation}
\limsup_{\rho\to\infty}\varepsilon(\rho)<2\,. \label{varepsilon}
\end{equation}
Then, given a sufficiently small $\delta>0$, we have
\[zE(z)+\frac1{\gamma (0)} = \left(1+o(1)\right)\sqrt{2\pi\frac{L(s)}{L'(s)}}
\exp\left(s^2\frac{L'(s)}{L(s)}\right)+o(1), \quad z\to\infty,\]
uniformly in $\Omega(\pi/2+\delta)$, and
\[zE(z)+\frac1{\gamma (0)} = o(1),\quad z\to\infty \]
uniformly  in $\mathbb{C}\setminus \Omega(\pi/2+\delta)$.
Here, also  $s=s_z$ and the   branch of the  square root is  positive on the positive half-line.
\end{theorem}
Note that it is not difficult to drop assumption~\eqref{varepsilon}
in Theorem~\ref{TheoremE} at the expense of a more complicated conclusion.
We shall not do this here. One of the reasons is that we are mainly interested in the case when $\varepsilon(\rho)=o(1)$ as $\rho \to\infty$.

\medskip
We also note that the asymptotics given in Theorems~\ref{TheoremE} and~\ref{TheoremK}
are known in the case when there exists a positive limit
\[ \lambda=\lim_{r\to\infty} \varepsilon(\rho), \]
cf. \cite{Evgrafov4, Mayergoiz}. In this case, $E$ is an entire function of order $\frac{1}{\lambda}$. The logarithmic case $L(s)=\log(s+e)$ is also classical and goes back to Lindel\"of.

\subsubsection{An example.} $L(s)=\log^\beta(s+e),$ $\beta>0$. In this case, the saddle-point equation \eqref{saddlePoint} has the form
\[ \beta\log\log (s+e)+\frac{\beta s}{s+e}\cdot\frac{1}{\log(s+e)}=\log z, \]
which readily simplifies to
\[\log\log s+\frac{1}{\log s}+O\left(\frac{1}{s\log s}\right)=\frac{1}{\beta}\log z,  \]
whence
\[ s=\left(1-\frac{1+o(1)}{2}z^{-1/\beta}\right)\exp\left(z^{1/\beta}-1\right). \]
Then
\[ K(z)=\left(1+o(1)\right)\sqrt{\frac{s\log s}{2\pi \beta}}\exp\left(-\beta\frac{s}{\log s}\right) \]
uniformly in any domain $\Omega(\pi-\delta)$, and
\[ zE(z)+1=\left(1+o(1)\right)\sqrt{\frac{2\pi}{ \beta}s\log s}\,\exp\left(\beta\frac{s}{\log s}\right) + o(1) \]
uniformly in $\Omega(\pi/2+\delta)$ with sufficiently small $\delta>0$.

\medskip
We note that the entire function $E$ has nearly maximal growth in the curvilinear strip $\Omega(\pi/2)$, while the analytic function $K$ has nearly fastest decay in $\Omega(\pi/2)$, and that, for sufficiently large $r_0$,
\[ \Omega(\pi/2)\cap \{z: |z|>r_0\}=\{z=re^{i\psi}:\;|\psi|\leq \Psi(r),\;r>r_0 \}, \]
where
\[ 
\Psi(r)=\frac{\pi\beta}{2}\left(r^{-1/\beta} + 
\left( \tfrac{\pi^2}8-\tfrac12 \right) r^{-3/\beta}+O(r^{-4/\beta})\right),\quad r\to\infty. \]

\subsubsection{} The observation we have just made is quite general, For any curvilinear semistrip $\Omega$ which is bounded by two sufficiently regularly varying curves $\{z=re^{\pm i \Psi(r)}\}$, one can find a function $\gamma$, satisfying our regularity conditions (A), (B) and (C), so that the entire function $E$ will have nearly maximal growth in $\Omega$, while the  analytic function $K$ will have nearly fastest decay in $\Omega$.  We shall not pursue that matter here.

\subsection{Examples of admissible functions $\gamma$.}
We start with several straightforward observations:

\subsubsection{} The shifted Euler's Gamma-function $\Gamma(s+c)$, $c>0$, is admissible.

\subsubsection{} If the function $\gamma$ is admissible, then the functions
\[
s\mapsto\frac{\gamma(s+c)}{\gamma(c)},\ c>0, \quad \text{and} \quad
s\mapsto \gamma(s)e^{\tau s}, \ \tau\in\mathbb{R}\,,
\]
are also admissible.

\subsubsection{} Denote by $\log_k$ the $k$-th iterate of the logarithmic function. Then the function
\[
\gamma(s)=\exp\left[a s\log_k^b(s+c_k)\right],\ k\in\mathbb{N}\,,
\]
is admissible provided that $a>0$, $b>0$ (and $b\leq 1$ for $k=1$), and that $c_k>0$ are sufficiently large.

\medskip
The following simple rules allow one to construct new admissible functions form the given ones:
\subsubsection{} If $\gamma$ is admissible and $a>0$, then the function $\gamma^a$ is also admissible.
\subsubsection{} If $\gamma_1$ and $\gamma_2$ are admissible, then $\gamma_1\cdot\gamma_2$ is always admissible, while $\frac{\gamma_1}{\gamma_2}$ is admissible provided that $\gamma_1\geq \gamma_2$ on $(0,\infty)$ and that the function
\[ \rho\mapsto\left(\frac{\gamma_1(\rho)}{\gamma_2(\rho)}\right)^{1/\rho},\quad \rho>0, \]
is non-decreasing and unbounded.

\subsubsection{} If $\gamma(s)=L(s)^s$ is admissible, then the function
\[ s\mapsto\left(\log L(s+1)\right)^s \] is admissible as well.

\subsection{Admissible functions with prescribed asymptotic behavior.} It is not difficult to construct admissible functions with prescribed asymptotic behavior on the positive ray.
The next result is a version of the known observation that if $h$ is a slowly varying
function on $[0, \infty)$ (that is, $\rho h'(\rho)=o(h(\rho))$ as $\rho\to\infty$), then the function
\[
\mathfrak h(s) = s \int_0^\infty \frac{h(u)\, {\rm d}u}{(u+s)^2}
\]
is analytic in $\{s\colon |\arg(s)|<\pi\}$, slowly varying on $\mathbb R_+$, and for $\rho\to\infty$ satisfies $\mathfrak h(\rho e^{{\rm i}\theta}) = (1+o(1))h(\rho)$ uniformly
in any angle $|\arg(s)|\le \pi-\delta$.

\begin{theorem}\label{thm3}
	Suppose $\ell:[0,\infty)\to(0,\infty)$ is an unboundedly increasing $C^1$-function such that the function
	\[ \rho\mapsto \rho\frac{\ell'(\rho)}{\ell(\rho)} \]
	is slowly varying and bounded for $\rho>0$. Then, for any $c>0$, the function
	\begin{equation}
	 \gamma(s)=\exp\left(s^2\int_c^\infty \frac{\ell'(u)}{\ell(u)} \frac{{\rm d}u}{s+u} \right),\quad |\arg(s+c)|<\pi\,, \label{Postivetype}
	\end{equation}
	is admissible and
	\[ \lim_{\rho\to\infty}\frac{\log \gamma(\rho)}{\rho\log \ell(\rho)}=1. \]
	If, in addition, there exists the limit, $\displaystyle\lim_{\rho\to\infty} \rho\frac{\ell'(\rho)}{\ell(\rho)}$, then
	\[ \lim_{\rho\to\infty}\frac{\ell(\rho)}{ \gamma(\rho)^{1/\rho}}=\ell(0). \]
\end{theorem}
For the reader's convenience, we give the proof of Theorem~3 in Appendix~A.

\subsection{Admissible functions of positive type.}

Beurling observed in \cite{Beurling2,Beurling} that analytic functions that admit integral representations similar to \eqref{Postivetype} have special  positivity properties which yield that $K(t)\geq 0$ for $t>0$. This provides us with a large class of explicit integral representations for solutions $K(t)$ to the Stieltjes moment problem with known asymptotics as $t\to\infty$ given by Theorem~\ref{TheoremK}. We shall discuss this in Appendix~B.

\subsubsection*{Acknowledgments.} We thank Andrei Iacob for his help with copy-editing of this
paper.

\section{Preliminaries}\label{sect:preliminaries}

Put
\[
G(z,s):=\log\gamma(s)-s\log z=s\log L(s)-s\log z\,.
\]
Then
\begin{equation}\label{eqK-G}
K(z) = \frac1{2\pi {\rm i}}\, \int_{\mathcal L_\alpha} e^{G(z, s)}\, {\rm d}s,
\end{equation}
where $\mathcal L_\alpha$ is the same contour as on Figure~1, and
\begin{equation}\label{eqE-G}
z E(z) + \frac1{\gamma (0)} = \int_{-\sigma_0}^\infty e^{-G(z, s)}\, {\rm d}s + o(1)\,.
\end{equation}
The latter relation easily follows from the classical Abel-Plana summation
formula (see Section~4). In this section, we collect estimates of the function $G$ and its derivatives needed for the asymptotic estimates of the integrals on the RHS of~\eqref{eqK-G}
and~\eqref{eqE-G}.

\medskip
Recall that, for any $0<\alpha<\alpha_0$, the function $\Phi (s) = \log L(s) + \varepsilon (s)$ (where, as before, $\varepsilon (s) = s\tfrac{L'}L(s)$) is univalent in the domain
$S(\alpha, \rho_0) = \bigl\{ s\colon |\arg (s)|<\alpha, |s|>\rho_0 \bigr\}$ with sufficiently large $\rho_0$,
and that we denote $\Omega(\alpha, \rho_0) = \Phi(S(\alpha, \rho_0))$. Hence, for any $z\in\Omega(\alpha, \rho_0)$, the function $s\mapsto G(z, s)$ has a unique critical point, which we
denote by $s_z=\rho_ze^{{\rm i}\theta_z}$, and this is the unique saddle-point of the function $s\mapsto \re G(z, s)$.

\subsection{Derivatives of $G$.}
The first derivative $G'_s(z, s)$ equals $\log L(s) + \e (s) - \log z$. Recalling conditions
(D) and (C), and the saddle-point equation
\begin{align*}
\log z &= \log L(s_z) + \e (s_z) \\ \\
&= \int_0^{\rho_z} \frac{\e (u)}{u}\, {\rm d}u + \e (\rho_z) + {\rm i}\theta_z \e (\rho_z)
+ o(\e (\rho_z))\,,
\end{align*}
we get
\begin{equation}\label{eq:G'}
G'_s(z, \rho e^{{\rm i}\theta}) = \int_{\rho_z}^\rho \frac{\e (u)}{u}\, {\rm d}u
+ {\rm i}\e (\rho_z) (\theta-\theta_z) + \chi_1(z) + \chi_2(\rho e^{{\rm i}\theta})\,,
\end{equation}
where $\chi_1 (z) = o(\e (\rho_z))$ uniformly in $z\in\bar\Omega(\alpha_0-\delta)$, $z\to\infty$, and $\chi_2(\rho e^{{\rm i}\theta}) = o(\e (\rho))$ uniformly in
$|\theta|\le \alpha_0-\delta$, $\rho\to\infty$.

\medskip
The second derivative $G''_{ss}(z, s)$ does not depend on $z$ and equals
\begin{align}\label{eq:G''}
G''_{ss}(z, s) &= \frac{L'}L(s) + \varepsilon'(s) \nonumber \\
&= (1+o(1)) \frac{\varepsilon(\rho)}s \qquad (\text{by}\ (B)\ \text{and}\ (C)\,)
\end{align}
uniformly in any angle $|\arg(s)|\le \alpha_0-\delta$. In particular,
\begin{equation}\label{eq:G''(s_z)}
G''_{ss}(z, s_z) = (1+o(1))\, \frac{\varepsilon(\rho_z)}{\rho_z}\, e^{-{\rm i}\theta_z}\,.
\end{equation}
Since the function $\varepsilon (\rho)$ is slowly varying, we see that if we will
succeed to correctly deform the integration contours in the integrals on the RHSs
of~\eqref{eqK-G} and~\eqref{eqE-G}, then the asymptotics of these integrals will
be determined by a neighbourhood of the saddle point $s_z$ of size $\rho_z^c$ with any $c>\tfrac12$.

\subsection{Behaviour of $G$ in a neighbourhood of the saddle point $s_z$.}

We fix a small positive $\delta_1<\tfrac12$ (for instance, the value $\delta_1=\tfrac18$
will suffice for our purposes) and assume that $|s-s_z|\le \rho_z^{1-\delta_1}$. By~\eqref{eq:G''},
combined with condition (B), we have
\begin{equation}\label{eq:G''-local}
G''_{ss}(z, s) = (1+o(1))\, G''_{ss}(z, s_z)
\end{equation}
uniformly in $|s-s_z|\le \rho_z^{1-\delta_1}$, whence,
\begin{equation}\label{eq:G-local}
G(z, s) = G(z, s_z) + \bigl( \tfrac12+o(1) \bigr)\, (s-s_z)^2\, \frac{\varepsilon (\rho_z)}{\rho_z}\, e^{-{\rm i}\theta_z}
\end{equation}
also uniformly in $|s-s_z|\le \rho_z^{1-\delta_1}$, $z\in\bar\Omega(\alpha_0-\delta)$, $z\to\infty$. Thus,
\begin{itemize}
\item
{\em the function $w\mapsto \re G(z, s_z + w)$ has the fastest decay in the directions $w=\pm {\rm i}e^{{\rm i}\theta_z/2}$}
\end{itemize}
and
\begin{itemize}
\item
{\em the fastest growth in the directions
$w=\pm e^{{\rm i}\theta_z/2}$}.
\end{itemize}

\medskip Let $\Gamma$ be a smooth simple curve that traverses once the disk
\[ D(s_z)=\{|s-s_z|\le \rho_z^{1-\delta_1} \} \]
and passes through the saddle-point $s_z$.
\begin{figure}[h]
	\includegraphics[scale=0.85]{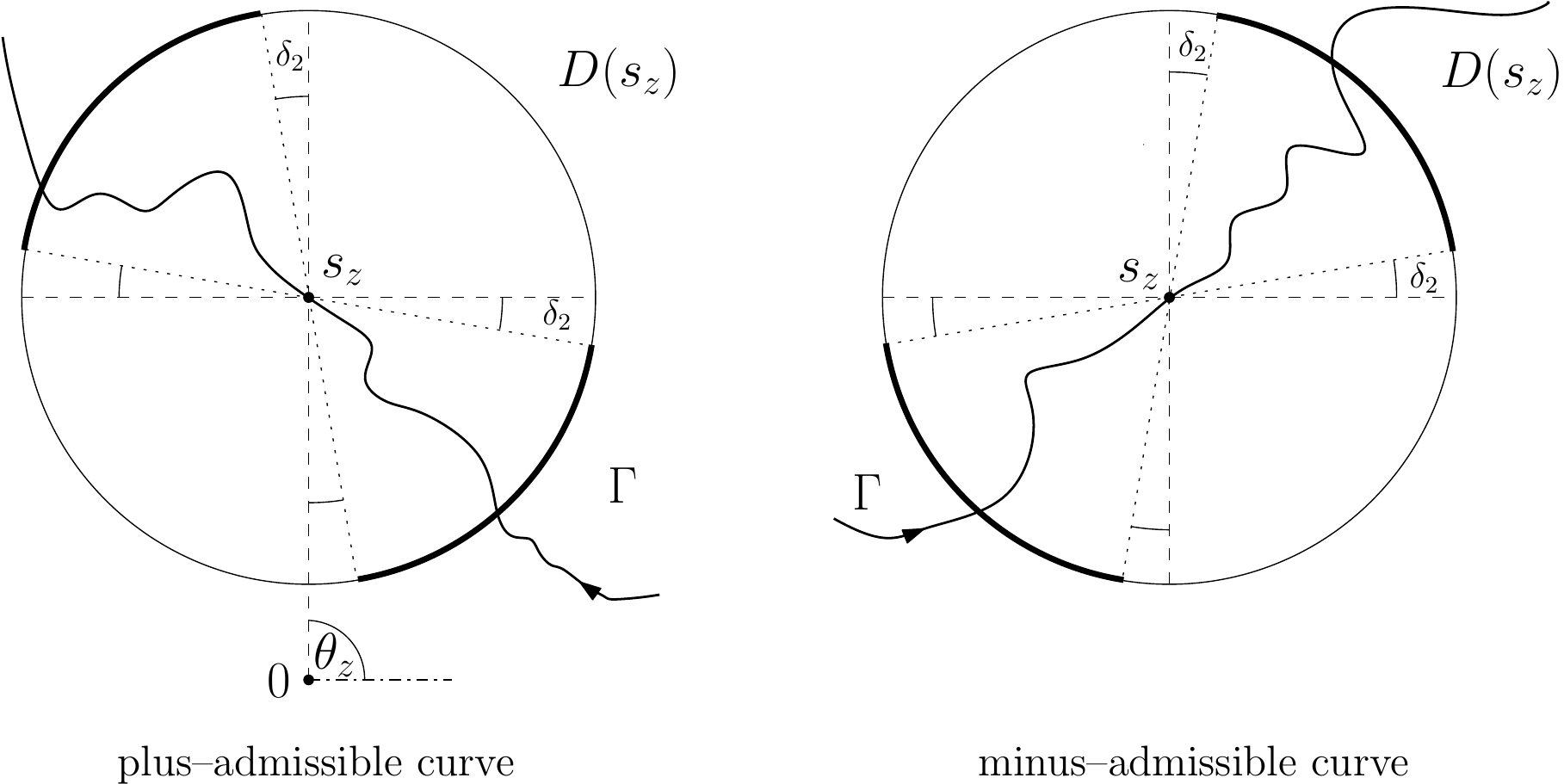}
\end{figure}
We call the curve $\Gamma$ {\em plus-admissible} if it
\begin{itemize}
\item enters $D(s_z)$ through the
arc $\{s=s_z-{\rm i} \rho_z^{1-\delta_1}e^{{\rm i}\varphi}\colon \bigl| \varphi - \tfrac12 \theta_z \bigr| \le \tfrac14\pi - \delta_2\}$,
\item exists $D(s_z)$ through the arc
$\{s=s_z+{\rm i} \rho_z^{1-\delta_1}e^{{\rm i}\varphi}\colon \bigl| \varphi - \tfrac12 \theta_z \bigr| \le \tfrac14\pi - \delta_2 \}$,
\end{itemize}
and
\begin{itemize}
\item
$\Gamma \bigcap D(s_z)$ does
not leave the set
\[
\bigl\{ s=s_z+{\rm i} t e^{{\rm i}\varphi}\colon
-\rho_z^{1-\delta_1}\le t \le \rho_z^{1-\delta_1},\, \bigl| \varphi - \tfrac12 \theta_z \bigr| \le \tfrac14\pi - \delta_2 \bigr\}\,.
\]
\end{itemize}

Similarly, we say that $\Gamma$ is {\em minus-admissible}, if
\begin{itemize}
\item
it enters the disk $D(s_z)$ through the arc
$\{ s=s_z- \rho_z^{1-\delta_1}e^{{\rm i}\varphi}\colon \bigl| \varphi - \tfrac12 \theta_z \bigr| \le \tfrac14\pi - \delta_2\}$,
\item
exits $D(s_z)$ through the arc
$\{s=s_z+\rho_z^{1-\delta_1}e^{{\rm i}\varphi}\colon \bigl| \varphi - \tfrac12 \theta_z \bigr| \le \tfrac14\pi - \delta_2\}$,
\end{itemize}
and
\begin{itemize}
\item
$\Gamma\bigcap D(s_z)$ does not leave the set
\[
\bigl\{ s=s_z + t e^{{\rm i}\varphi}\colon
-\rho_z^{1-\delta_1}\le t \le \rho_z^{1-\delta_1},\, \bigl| \varphi - \tfrac12 \theta_z \bigr| \le \tfrac14\pi - \delta_2 \bigr\}\,.
\]
\end{itemize}

\begin{lemma}\label{lemma-s.p.}
Suppose that the curve $\Gamma$ is plus-admissible. Then
\[
\int_{\Gamma\cap D(s_z)} e^{G(z, s)}\, {\rm d}s = ({\rm i}+o(1))\, \sqrt{2\pi\, \tfrac{L}{L'}(s_z)}\,
e^{-s_z^2\, \frac{L'}{L}(s_z)}\,.
\]
If the curve $\Gamma$ is minus-admissible, then
\[
\int_{\Gamma\cap D(s_z)} e^{-G(z, s)}\, {\rm d}s = (1+o(1))\, \sqrt{2\pi\, \tfrac{L}{L'}(s_z)}\,
e^{s_z^2\, \frac{L'}{L}(s_z)}\,.
\]
Both asymptotic relations hold
uniformly in $z\in\bar\Omega(\alpha_0-\delta)$, $z\to\infty$. The branch of the square root on the RHSs is positive when $z$ belongs to the positive ray.
\end{lemma}

\noindent
Note that $s_z^2 L'(s_z)/L(s_z) = s_z \e (s_z)$, and $L(s_z)/L'(s_z)
= s_z/\e (s_z) = (1+o(1))s_z/\e (\rho_z)$.

\medskip\noindent{\em Proof of Lemma~\ref{lemma-s.p.}}:
We prove only the first statement; the
proof of the second one is very similar.
We start with the special case, when $\Gamma\cap D(s_z)$ is the segment $I$ of the
fastest decay of the function $w\mapsto \re G(z, s_z+w)$:
\[
I = \bigl\{ s=s_z+{\rm i}te^{{\rm i}\theta_z/2}\colon -\rho_z^{1-\delta_1} \le t \le \rho_z^{1-\delta_1} \bigr\}\,.
\]
In this case,
\begin{align*}
\int_{I} e^{G(z, s)}\, {\rm d}s &\stackrel{\eqref{eq:G-local}}= e^{G(z, s_z)}\, \int_{I_+} e^{(\frac12+o(1))(s-s_z)^2\,
\frac{\varepsilon (\rho_z)}{\rho_z}e^{-{\rm i}\theta_z}}\, {\rm d}s \\ \\
&= e^{G(z, s_z)}\, {\rm i}e^{{\rm i}\theta_z/2}\, \sqrt{\frac{\rho_z}{\varepsilon(\rho_z)}}\,
\int_{-\rho_z^{1/2-\delta_1}\varepsilon(\rho_z)^{-1/2}}^{\rho_z^{1/2-\delta_1}\varepsilon(\rho_z)^{-1/2}} e^{-(\frac12+o(1))t^2}\, {\rm d}t\,.
\end{align*}
The function $\varepsilon (\rho)$ is slowly varying. Hence, for any $c>0$,
the function $\rho^{-c}\varepsilon(\rho)$ decays to $0$ as $\rho\to\infty$.
Since $\delta_1<\tfrac12$, we conclude that
$ \rho_z^{1/2-\delta_1}\varepsilon(\rho_z)^{-1/2} \to +\infty $ uniformly in
$z\in\bar\Omega(\alpha_0-\delta)$, $z\to\infty$, and therefore, the integral on the
RHS converges to $\sqrt{2\pi}$ also uniformly.

It remains to note that $G(z, s_z)=-s_z^2 \tfrac{L'}L(s_z)$, and that, by condition (C),
\[
e^{-{\rm i}\theta_z/2}\, \sqrt{\frac{\rho_z}{\varepsilon(\rho_z)}}
=(1+o(1))\sqrt{\frac{L}{L'}(s_z)}\,,
\]
completing the proof of this special case of Lemma~\ref{lemma-s.p.}.

\medskip To move to the general case, we note that on the boundary
circumference $|s-s_z|=\rho_z^{1-\delta_1}$ the function $\re G(z, s)$ is
much smaller than at the saddle-point $s_z$.
More precisely, we
claim that {\em given small positive $\delta$, $\delta_1$, $\delta_2$, there exists a sufficiently large
$\rho_0=\rho_0(\delta, \delta_1, \delta_2)$ so that,
for $z\in\bar\Omega(\alpha_0-\delta, \rho_0)$ and $\bigl| \varphi-\tfrac12(\theta_z\pm \pi) \bigr| \le \tfrac14\pi -\delta_2$, we have}
\[
\re G\bigl( z, s_z+\rho_z^{1-\delta_1}e^{{\rm i}\varphi} \bigr)
\le \re G(z, s_z) - \rho_z^{1-3\delta_1}\,.
\]
Indeed, for $\bigl| \varphi-\tfrac12(\theta_z\pm \pi) \bigr| \le \tfrac14\pi -\delta_2$,
we have $\cos (2\varphi-\theta_z)\le - c <0$, whence, taking into account~\eqref{eq:G-local},
\begin{align*}
\re G(z, s) &\le \re G(z, s_z) - \bigl( \frac12 + o(1) \bigr)\, \rho_z^{2-2\delta_1}\, \frac{\varepsilon (\rho_z)}{\rho_z} \cdot c \\ \\
&\le \re G(z, s_z) - c_1 \rho_z^{1-2\delta_1} \varepsilon (\rho_z)\,,
\end{align*}
provided that $\rho_z$ is sufficiently large. Since the function $\varepsilon (\rho)$ is slowly varying, for $\rho$ sufficiently large,
we have $c_1\cdot \varepsilon (\rho) \ge \rho^{-\delta_1}$, proving the claim.

\medskip At last, $\rho_z^{1-3\delta_1} e^{-\rho_z^{1-3\delta_1}}$ is much smaller than $\displaystyle \sqrt{\frac{\rho_z}{\varepsilon (\rho_z)}}$. Thus,
using Cauchy's theorem, we can replace the segment $I$ by $\Gamma\bigcap D(s_z)$ for any plus-admissible curve $\Gamma$, completing the proof. \hfill $\Box$

\subsection{Asymptotics of $\re G$.}

We have
\begin{align*}
\re G(re^{{\rm i}\psi}, \rho e^{{\rm i}\theta}) &= \re \bigl[ \rho e^{{\rm i}\theta}
\log L(\rho e^{{\rm i}\theta}) - \rho e^{{\rm i}\theta} (\log r + {\rm i}\psi ) \bigr] \\ \\
&\stackrel{(D)}=
\rho\cos\theta \Bigl( \int_0^\rho \frac{\e (u)}u\, {\rm d}u - \log r \Bigr)
- \rho\sin\theta (\theta\varepsilon (\rho) - \psi) + o(\rho\varepsilon (\rho))
\end{align*}
uniformly in $|\theta|\le \alpha_0-\delta$, $\rho\to\infty$. Recalling the equation
\[
\log z = \log L(s_z) +\varepsilon (s_z)
\]
for the saddle point $s_z=\rho_ze^{{\rm i}\theta_z}$ and using conditions (C) and (D),
we see that
\begin{align}\label{eq:ReG-asymptotics}
\nonumber
\re G(z, \rho e^{{\rm i}\theta}) &= \rho\cos\theta \Bigl( \int_{\rho_z}^\rho
\frac{\e (u)}u\, {\rm d}u - \e(\rho_z) \Bigr) \\ \nonumber \\
&\qquad - \rho\sin\theta \bigl( \theta\e(\rho) - \theta_z\e(\rho_z) \bigr)
+\rho \bigl( \chi_1 (z) + \chi_2 (\rho e^{{\rm i}\theta}) \bigr)\,,
\end{align}
where, as before,
$\chi_1 (z) = o(\e (\rho_z))$ uniformly in $z\in\bar\Omega(\alpha_0-\delta)$, $z\to\infty$, and $\chi_2(\rho e^{{\rm i}\theta}) = o(\e (\rho))$ uniformly in
$|\theta|\le \alpha_0-\delta$, $\rho\to\infty$.

\subsection{Estimates of $\re G$ on arcs of the circumference $|s|=\rho_z$.}

\begin{lemma}\label{lemma:ReG on arcs}
Suppose that $ \max(|\theta|, |\theta_z| )\le \alpha_0-\delta$.
Then
\begin{equation}\label{eq:ReG-arcs1}
\re G(z, \rho_z e^{{\rm i}\theta}) - \re G(z, \rho_z e^{{\rm i}\theta_z})
= -(f(\theta, \theta_z)) + o(1)) \rho_z \e (\rho_z)\,,
\end{equation}
were $f(\theta, \theta_z) = \cos\theta - \cos\theta_z  + (\theta-\theta_z)\sin\theta$, and
\begin{equation}\label{eq:ReG-arcs2}
\frac{\partial^2}{\partial\theta^2}\, \re G(z, \rho_z e^{{\rm i}\theta})
 = (h(\theta) + o(1)) \rho_z \e (\rho_z)\,,
\end{equation}
were $h(\theta) = - \cos\theta + (\theta-\theta_z)\sin\theta$. Both estimates are uniform
as $z\to\infty$.
\end{lemma}

\noindent{\em Proof of Lemma~\ref{lemma:ReG on arcs}}: Estimate~\eqref{eq:ReG-arcs1} immediately follows from asymptotics~\eqref{eq:ReG-asymptotics}. The proof of~\eqref{eq:ReG-arcs2} is also straightforward. We have
\[
\frac{\partial^2}{\partial\theta^2}\, G(z, \rho_z e^{{\rm i}\theta})
= - \rho_z e^{{\rm i}\theta}\, G'_s(z, \rho_z e^{{\rm i}\theta})
- \rho_z^2 e^{2{\rm i}\theta}\, G''_{ss}(z, \rho_z e^{{\rm i}\theta})\,,
\]
\[
G'_s(z, \rho_z e^{{\rm i}\theta})
\stackrel{\eqref{eq:G'}}= {\rm i}(\theta-\theta_z)\e (\rho_z) + o(\e (\rho_z)\,,
\]
and
\[
G''_{ss}(z, \rho_z e^{{\rm i}\theta})
\stackrel{\eqref{eq:G''}}= (1+o(1)) \frac{\e (\rho_z)}{\rho_z}\, e^{-{\rm i}\theta_z}\,.
\]
Therefore,
\begin{align*}
\frac{\partial^2}{\partial\theta^2}\, \re G(z, \rho_z e^{{\rm i}\theta}) &=
\re \Bigl[\frac{\partial^2}{\partial\theta^2}\, G(z, \rho_z e^{{\rm i}\theta})\Bigr] \\ \\
&=
- \rho_z \e(\rho_z)\,
\re\Bigl[ 1+ e^{{\rm i}\theta} \bigl( {\rm i}(\theta-\theta_z) + o(1)\bigr) \Bigr] \\ \\
&= - \rho_z \e(\rho_z)\, \bigl( \cos\theta - (\theta-\theta_z) \sin\theta + o(1)\bigr)\,,
\end{align*}
completing the proof. \hfill $\Box$

\subsection{Estimate of $\re G$ on segments that pass through the saddle point.}

\begin{lemma}\label{lemma:ReG-segment}
Let $t$ be a real number such that $|t|\le 1-\delta_3$. Then,
\[
\frac{\partial^2}{\partial t^2}\, \re G(z, s_z+t\rho_z e^{{\rm i}\ph})
= (1+o(1)) \frac{\rho_z \e(\rho_z)}{|e^{{\rm i}\theta_z} + te^{{\rm i}\ph}|}\,
\cos\bigl( 2\ph - \arg(e^{{\rm i}\theta_z} + te^{{\rm i}\ph})\bigr)\,.
\]
\end{lemma}

\noindent{\em Proof of Lemma~\ref{lemma:ReG-segment}}:
\begin{align*}
\frac{\partial^2}{\partial t^2}\, \re G(z, s_z+t\rho_z e^{{\rm i}\ph})
&= \re \Bigl[ \rho_z^2 e^{2{\rm i}\ph} G''_{ss} (z, s_z+t\rho_z e^{{\rm i}\ph}) \Bigr] \\ \\
&\stackrel{\eqref{eq:G''}}=
\re \Bigl[ \rho_z^2 e^{2{\rm i}\ph} \cdot (1+o(1))\,
\frac{\e (s_z+t\rho_z e^{{\rm i}\ph})}{s_z+t\rho_z e^{{\rm i}\ph}}\Bigr]
\end{align*}
By (C) and (B), the RHS equals
\[
(1+o(1)) \frac{\rho_z \e(\rho_z)}{|e^{{\rm i}\theta_z} + te^{{\rm i}\ph}|}\,
\cos\bigl( 2\ph - \arg(e^{{\rm i}\theta_z} + te^{{\rm i}\ph})\bigr)\,
\]
which proves the lemma. \hfill $\Box$

\subsection{Tail estimates of $\re G$.}

\begin{lemma}\label{lem:ReG-tails}
The function
$\displaystyle \rho\mapsto \frac{\partial\re G(z, \rho e^{{\rm i}\theta})}{\partial \rho} $,
$\rho\ge\rho_z$, increases whenever $|\theta|\le \frac{\pi}2 - \delta_4$, and decays
whenever $|\theta|\ge \frac{\pi}2 + \delta_4$.
Furthermore, for $\rho_z\le \rho \le 2\rho_z$,
\[
\frac{\partial\re G(z, \rho e^{{\rm i}\theta})}{\partial \rho} =
\e (\rho_z) \Bigl( \log \frac{L(\rho)}{L(\rho_z)} \cos\theta- (\theta-\theta_z)\sin\theta
+o(1) \Bigr)
\]
uniformly in $z\in\bar\Omega(\alpha_0-\delta)$, $z\to\infty$.
\end{lemma}

\noindent{\em Proof of Lemma~\ref{lem:ReG-tails}}:
The first claim follows since, due to the asymptotics~\eqref{eq:G''},
we have
\[
\frac{\partial^2 \re G(z, \rho e^{{\rm i}\theta})}{\partial \rho^2}
= (1+o(1))\, \frac{\e (\rho)}{\rho}\, \cos\theta\,.
\]

To check the second claim we write
\begin{align*}
\frac{\partial\re G(z, \rho e^{{\rm i}\theta})}{\partial \rho} &=
\re \Bigl[ e^{{\rm i}\theta} G'_s(z, \rho e^{{\rm i}\theta}) \Bigr] \\ \\
&\stackrel{\eqref{eq:G'}}= \re \Bigl[ e^{{\rm i}\theta}
\Bigl( \int_{\rho_z}^\rho \frac{\e (u)}u\, {\rm d}u + {\rm i} \e (\rho_z) (\theta-\theta_z) +o(\e(\rho_z)) \Bigr) \Bigr] \quad (\text{by\ } (C) \text{\ and\ } (D)) \\ \\
&= \e (\rho_z) \re \Bigl[ e^{{\rm i}\theta}
\Bigl( \log\frac{L(\rho)}{L(\rho_z)} + {\rm i} (\theta-\theta_z) +o(1) \Bigr) \Bigr] \\ \\
&= \e (\rho_z) \Bigl( \log \frac{L(\rho)}{L(\rho_z)} \cos\theta- (\theta-\theta_z)\sin\theta
+ o(1) \Bigr) \qquad \qquad (\text{by\ } (C))\,,
\end{align*}
proving the second claim. \hfill $\Box$

\subsection{Estimate of $\re G$ in a bounded sector.}

\begin{lemma}\label{lemma:reG-bounded-sector}
Given positive $\rho_1$ and $\delta$, there exists a positive $C$ so that
\[
\max \bigl\{ \bigl| \re G(z, \rho e^{{\rm i}\theta})
+ \re (s)\, \int_0^{\rho_z} \frac{\e (u)}{u}\, {\rm d} u
\bigr|\colon
\rho\le \rho_1, |\theta|\le \alpha-\delta \bigr\}
\le C\,,
\]
uniformly in $z\in\bar\Omega(\alpha_0-\delta)$.
\end{lemma}

\noindent{\em Proof of Lemma~\ref{lemma:reG-bounded-sector}}:
We have
\[
\re G(z, s) = \re (s) \cdot \log r + O(1)\,,
\]
uniformly in $s$, $|s|\le \rho_1$, $|\arg (s)|\le \alpha_0  - \delta$, and
in $z\in\bar\Omega(\alpha_0-\delta)$. Recalling that, by the saddle-point equation,
\[
\log r = \int_0^{\rho_z} \frac{\e (u)}{u}\, {\rm d}u + O(1)\,,
\]
also uniformly in $z\in\bar\Omega(\alpha_0-\delta)$, we get the result.
\hfill $\Box$

\section{Proof of Theorem~\ref{TheoremK}}

We fix several sufficiently small positive parameters $\delta$, $\delta_i$,
some of which already appeared in lemmas proven in the previous section.
Some restrictions on these parameters will be imposed in the course of the proof.
By $c$ and $C$ we denote various positive constants that may depend on these
parameters, the values of these constants are inessential for our purposes and
may differ from line to line.

In the course of the proof, all expressions that are
\[
o(1)\, \sqrt{\frac{\rho_z}{\e (\rho_z)}}\, e^{\re G(z, s_z)}
\]
will be called {\em negligible}.

\bigskip
Without loss of generality, we assume during the proof that the saddle
point $s_z$ lies in the sector $0\le \theta_z \le \alpha_0-\delta$, and split the proof
into three cases:

\bigskip\noindent
(I) $0\le \theta_z \le \tfrac{\pi}2-\delta_5$,

\bigskip\noindent
(II)
$\tfrac{\pi}2-\delta_5\le \theta_z \le \tfrac{\pi}2+\delta_5$,

\bigskip\noindent
and

\bigskip\noindent
(III) $ \tfrac{\pi}2+\delta_5 \le \theta_z \le \alpha_0-\delta$,

\bigskip\noindent
where $\delta, \delta_5 < \tfrac13 (\alpha_0-\tfrac{\pi}2)$.
In each of these three cases, using the asymptotics~\eqref{eq:ReG-asymptotics} of $\re G$,
we deform the original integration contour $\mathcal L_\alpha$ into a plus-admissible
contour $\Gamma_z$ that passes through the saddle point $s_z$. Then, Lemma~\ref{lemma-s.p.}
gives us the asymptotics of $\Gamma\bigcap D(s_z)$ which is always the main term, and
in each of the three cases we will need to show that the integral
over $\Gamma_z\setminus D(s_z)$ is negligible.

\subsection{Case~I: $0\leq\theta_z<\frac{\pi}{2}-\delta_5$.} We introduce the curve
\[
\Gamma_z = - J_1 + J_2 + J_3\,,
\]
where
\begin{align*}
J_1 &= \bigl\{ s=\rho e^{-{\rm i}\theta_0}\colon \rho \ge \rho_z \bigr\}, \\ \\
J_2 &= \bigl\{ s=\rho_z e^{{\rm i}\theta}\colon -\theta_0 \le \theta \le \theta_0 \bigr\}, \\ \\
J_3 &= \bigl\{ s=\rho e^{{\rm i}\theta_0}\colon \rho  \ge \rho_z \bigr\},
\end{align*}
with $\theta_0 = \tfrac{\pi}2+\delta_4$, and then split the arc $J_2$ into three parts,
\[ J_2 = J_+ \bigcup J_2' \bigcup J_2'',\]
where
\begin{align*}
J_+ &= \bigl\{s\in J_2\colon |s-s_z|\le \rho_z^{1-\delta_1} \bigr\}, \\ \\
J_2' &= \bigl\{s\in J_2\colon  |s-s_z|\ge \rho_z^{1-\delta_1}, |\theta-\theta_z|\le \delta_6 \bigr\}, \\ \\
J_2'' &= J_2 \setminus \bigl( J_+ \bigcup J_2' \bigr).
\end{align*}
\begin{figure}[h]
 \includegraphics[scale=0.75]{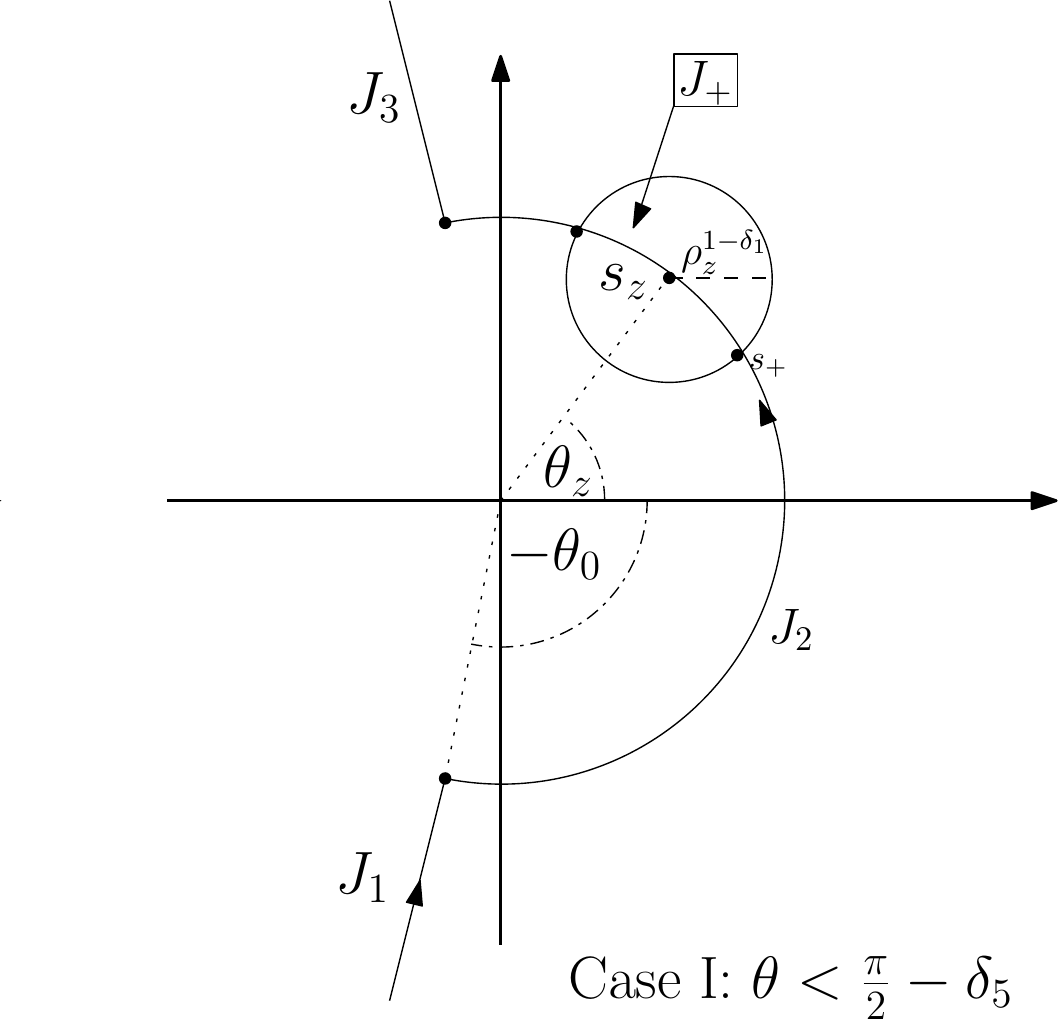}
   \caption{$\Gamma_{z} $}
\end{figure}
It is easy to see that the arc $J_2$ is plus-admissible, so
it remains to show that the integrals over $J_2'$, $J_2''$, $J_1$, and $J_3$ are
negligible.

\subsubsection{Integral over $J_2'$.}

Here, we will use estimate~\eqref{eq:ReG-arcs2} in Lemma~\ref{lemma:ReG on arcs}.
For $|\theta-\theta_z|\le \delta_6$, we have $\delta_6 \le \theta \le \tfrac{\pi}2 - \delta_5 + \delta_6$. Therefore,
\begin{multline*}
h(\theta) = - \cos\theta + (\theta-\theta_z)\sin\theta
\le -\cos \bigl( \frac{\pi}2 - \delta_5 + \delta_6 \bigr) + \delta_6 \\
\le - \frac{\delta_5-\delta_6}{\pi} + \delta_6 \le -c <0\,,
\end{multline*}
provided that $\delta_6 \le \tfrac12 \delta_5$. Hence, in this range,
\[
\frac{\partial^2}{\partial\theta^2}\, \re G(z, \rho_z e^{{\rm i}\theta})
\le -c \e(\rho_z)\rho_z\,,
\]
and then
\[
\re G(z, \rho_z e^{{\rm i}\theta}) \le \re G(z, s_z) - c\e (\rho_z)\rho_z
(\theta-\theta_z)^2\,.
\]
Thus,
\begin{align*}
\Bigl| \int_{J_2'} e^{G(z, s)}\, {\rm d}s \Bigr|
&\le e^{\re G(z, s_z)}\, \rho_z \int_{c\rho^{-\delta_1}_z\le|\theta-\theta_z|\le \delta_6}
e^{-c\e(\rho_z)\rho(z)(\theta-\theta_z)^2}\, {\rm d}\theta  \\ \\
&\le C e^{\re G(z, s_z)}\,
\sqrt{\frac{\rho_z}{\e(\rho_z)}\,} \cdot e^{-c(\e (\rho_z)\rho_z^{1-\delta_1})^2}
\end{align*}
with negligible RHS.

\subsubsection{Integral over $J_2''$.}\label{subsubsect:J_2''}

Now we will use estimate~\eqref{eq:ReG-arcs1} in Lemma~\ref{lemma:ReG on arcs}.
We claim that the function $f(\theta, \theta_z) = \cos\theta-\cos\theta_z + (\theta-\theta_z)\sin\theta$,
which appears on the RHS of~\eqref{eq:ReG-arcs1}
is strictly positive whenever
$|\theta|\le \tfrac{\pi}2+\delta_4$, $0\le \theta_z \le \tfrac{\pi}2-\delta_5$.
Indeed, since $f'_\theta (\theta, \theta_z)=(\theta-\theta_z)\cos\theta$,
the function $\theta\mapsto f(\theta, \theta_z)$ has a zero local minimum at
$\theta=\theta_z$ and two positive local maxima at $\theta=\pm \tfrac{\pi}2$.
Hence, it suffices to check that the values $f(\tfrac{\pi}2+\delta_4, \theta_z)$ and
$f(-\tfrac{\pi}2-\delta_4, \theta_z)$ are positive.
Since $f'_{\theta_z}(\theta, \theta_z) = -2\sin\theta$ is negative at
$\theta=\tfrac{\pi}2+\delta_4$ and positive at $\theta=-\tfrac{\pi}2-\delta_4$, we see that
$f(\tfrac{\pi}2+\delta_4, \theta_z) \ge f(\tfrac{\pi}2+\delta_4, \tfrac{\pi}2-\delta_5)$
and $f(-\tfrac{\pi}2-\delta_4, \theta_z) \ge f(-\tfrac{\pi}2-\delta_4, 0)$.
Finally, expanding in $\delta_4$ and $\delta_5$, we get
\begin{align*}
f(\tfrac{\pi}2+\delta_4, \tfrac{\pi}2-\delta_5)
&= -\sin\delta_4 - \sin\delta_5 + (\delta_4+\delta_5)\cos\delta_4 \\ \\
&= \tfrac16 (\delta_4^3 +\delta_5^3) - \tfrac12 \delta_4^2 (\delta_4+\delta_5)
+ O\bigl( \delta_4^4 + \delta_5^4 \bigr) \\ \\
&= \tfrac12 \bigl( \delta_4 + \delta_5 \bigr)\cdot
\bigl(\delta_5^2 - \delta_4\delta_5 + \tfrac12 \delta_4^2 \bigr)
+ O\bigl( \delta_4^4 + \delta_5^4 \bigr)
> 0\,,
\end{align*}
$\delta_4\le\tfrac12 \delta_5$, and
\begin{align*}
f(-\tfrac{\pi}2-\delta_4, 0) &= -\sin\delta_4 -1 + \bigl( \frac{\pi}2+\delta_4\bigr)\cos\delta_4 \\ \\
&=\frac{\pi}2 - 1 + O(\delta_4) >0\,.
\end{align*}
This proves the claim, which immediately
yields that on $J_2''$ we have $\re G(z, s) \le G(z, s_z) - c\rho_z\e(\rho_z)$, and therefore,
the integral over $J_2''$ is negligible.

\subsubsection{Integrals over $J_1$ and $J_3$.}\label{subsubsect:J_1-J_3}

Suppose that $s\in J_3$, that is, $s=\rho e^{{\rm i}\theta_0}$, $\rho\ge\rho_z$.
Then, by Lemma~\ref{lem:ReG-tails},
\[
\re G(z, \rho e^{{\rm i}\theta_0}) \le \re G(z, \rho_z e^{{\rm i}\theta_0})
- c \e(\rho_z)(\rho-\rho_z)\,.
\]
Besides, we already know that
\[
\re G(z, \rho_z e^{{\rm i}\theta_0}) \le  \re G(z, \rho_z e^{{\rm i}\theta_z})
- c_1 \e(\rho_z)\rho_z\,.
\]
Thus,
\[
\re G(z, \rho e^{{\rm i}\theta_0}) \le
\re G(z, \rho_z e^{{\rm i}\theta_z}) - c_1 \e(\rho_z)\rho_z
- c \e(\rho_z)(\rho-\rho_z)\,,
\]
and therefore, the integral over $J_3$ is negligible. For the same reason,
the integral over $J_1$ is negligible as well.

\subsection{Case II: $\frac{\pi}2-\delta_5 \le \theta_z \le \frac{\pi}2+\delta_5$.}
\label{subsect:KII}
We put, as in the previous case, $\theta_0 = \tfrac{\pi}2+\delta_4$ with
$\delta_4\ge 2\delta_5$. Consider the straight line
$\bigl\{ s=s_z+te^{{\rm i} 3\pi/4}\colon t\in\mathbb R\bigr\}$ and denote
by $s_+$ its intersection point with the ray $\{\arg (s)=\theta_0\}$ and
by $s_-\ne s_z$ its intersection point with the circumference $|s|=\rho_z$.
Let $J$ be the segment $[s_-, s_+]$. Clearly, it is a plus-admissible curve.
Our integration contour will be
\[
\Gamma_z = - J_1 + J_2 + J + J_3,
\]
where
\begin{align*}
J_1 &= \{s=\rho e^{-{\rm i}\theta_0}\colon \rho\ge \rho_z\}, \\ \\
J_2 &= \{s=\rho_z e^{{\rm i}\theta}\colon -\theta_0 \le \theta \le \arg (s_-)\}, \\ \\
J_3 &= \{s=\rho e^{{\rm i}\theta_0}\colon \rho \ge |s_+| \}.
\end{align*}
\begin{figure}[th]
\includegraphics[scale=0.75]{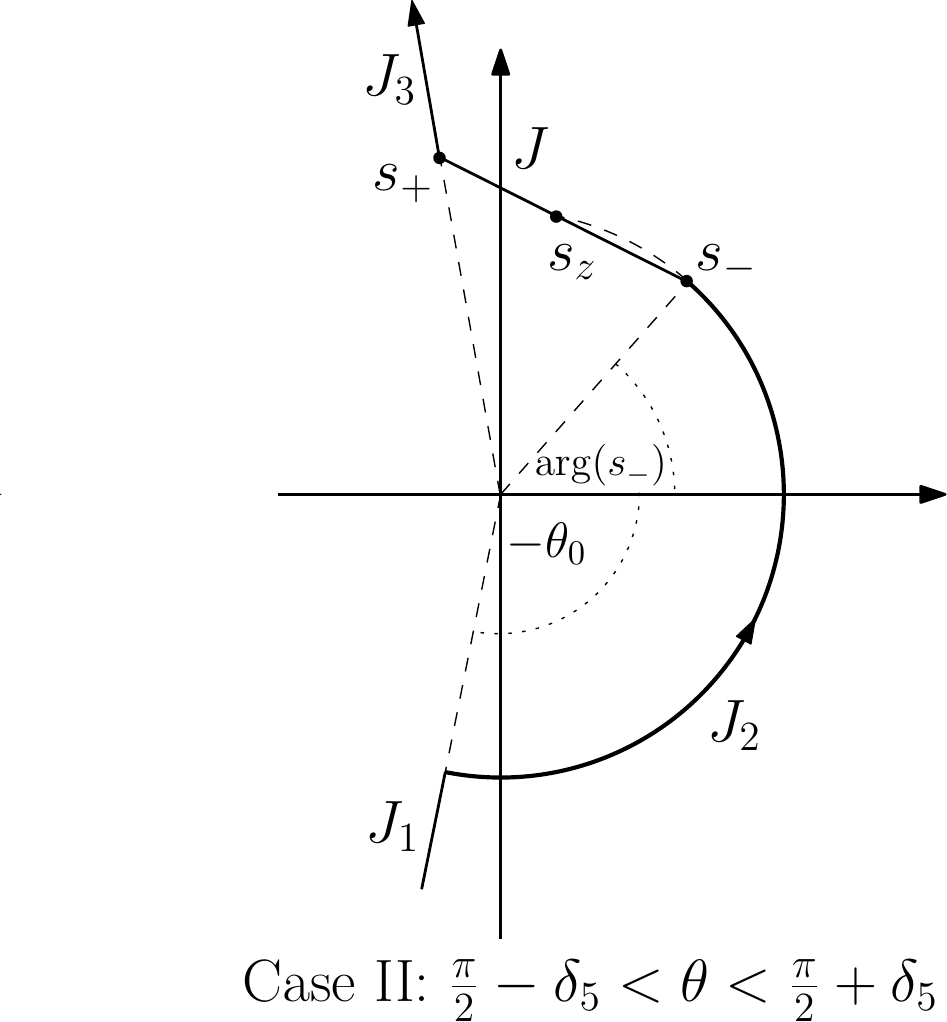}
   \caption{$\Gamma_{z} $}
\end{figure}

The main term in the asymptotics comes from the segment
$J_+=J\bigcap D(s_z)$, and we need
to check that the four remaining integrals over $J_1$, $J_2$, $J\setminus J_+$,
and $J_3$ are negligible. Estimates of the integrals over $J_1$ and $J_3$
follow the same lines as in~\ref{subsubsect:J_1-J_3}.
So here we estimate only the integrals over $J\setminus J_+$ and $J_2$.

\subsubsection{Integral over $J\setminus J_+$.}\label{subsubsect:KII-segment}
Here, we use Lemma~\ref{lemma:ReG-segment} with $\ph=\tfrac{3\pi}4$. Since,
for $s\in [s_-, s_+]$, $\cos(\tfrac{3\pi}2-\arg (s)) = -\sin(\arg s) \le -c<0$,
Lemma~\ref{lemma:ReG-segment} yields
\[
\frac{\partial^2}{\partial t^2}\, \re G(z, s_z+te^{{\rm i}\frac{3\pi}4}) \le
- c \rho_z \e(\rho_z)
\]
whenever $ s= s_z+te^{{\rm i}\frac{3\pi}4}\in J $,
whence,
\[
\re G(z, s_z+te^{{\rm i}\frac{3\pi}4}) \le \re G(z, s_z) - c\rho\e(\rho_z) t^2\,.
\]
Since we are on $J\setminus J_+$, we integrate only over
$|t|\ge c\rho_z^{1-\delta_1}$ and see that the integral over $J\setminus J_+$ is negligible.

\subsubsection{Integral over $J_2$.}

By estimate~\eqref{eq:ReG-arcs1}, we have
\[
\re G(z, \rho_z e^{{\rm i}\theta}) - \re G(z, \rho_z e^{{\rm i}\theta_z})
= - (f(\theta, \theta_z)+o(1))\rho_z\e(\rho_z)\,,
\]
with $f(\theta, \theta_z)=\cos\theta - \cos\theta_z +(\theta-\theta_z)\sin\theta$.
The same elementary analysis as in~\ref{subsubsect:J_2''}
shows that, for $ -\theta_0 \le \theta \le \theta_- =\arg(s_-)$ and
$\tfrac{\pi}2-\delta_5 \le \theta_z \le \tfrac{\pi}2+\delta_5$, one has
$f(\theta, \theta_z) \ge f(\theta_0, \theta_z)\ge c>0$. This implies
that the integral over the arc $J_2$ is negligible.

\subsection{Case~III: $\frac{\pi}2+\delta_5 \le \theta_z \le \alpha_0-\delta$.}
Here, we take the contour
\[
\Gamma_z = e^{{\rm i}\theta_z} \mathbb R_+ - e^{-{\rm i}\theta_z} \mathbb R_+ = J_1 - J_2\,.
\]
\begin{figure}[th]
\includegraphics[scale=0.75]{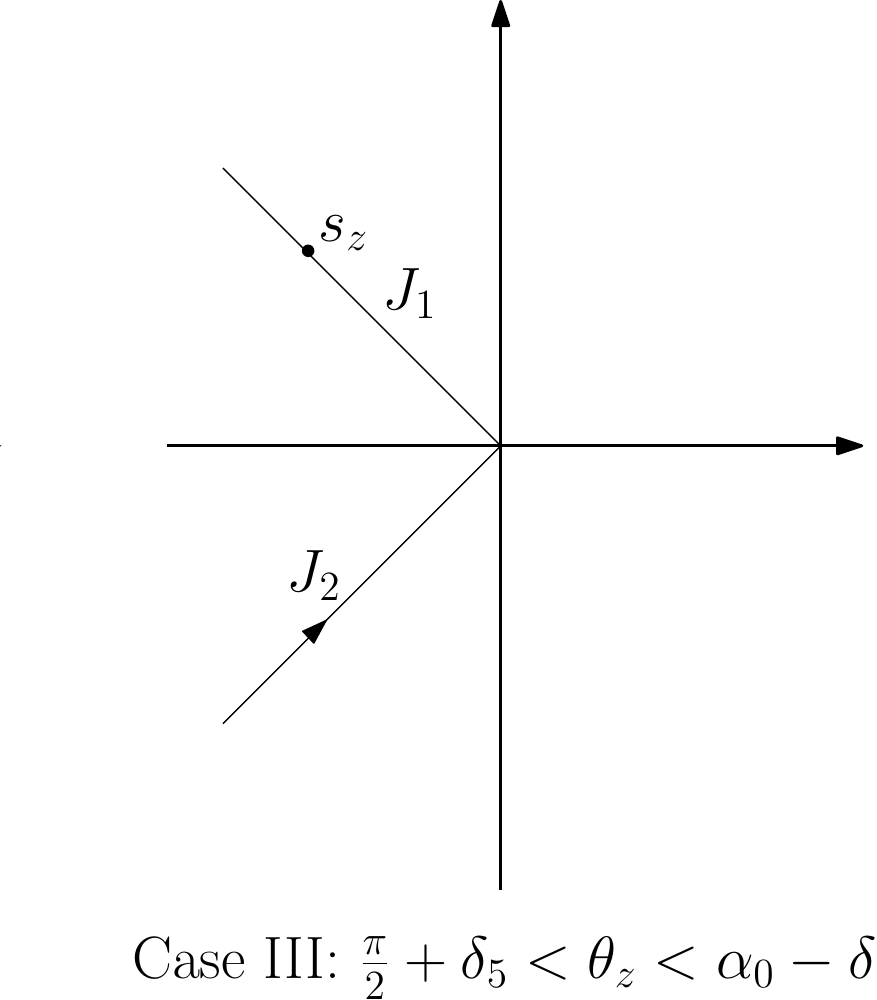}
   \caption{$\Gamma_{z} $}
\end{figure}
The ray $J_1=e^{{\rm i}\theta_z} \mathbb R_+$ is plus-admissible and
the main term in the asymptotics of the integral comes from integration
over the segment $J_1\bigcap D(s_z)$. Thus, we need to show that the integrals
over $J_1 \setminus D(s_z)$ and $J_2=e^{-{\rm i}\theta_z} \mathbb R_+$
are negligible.

\subsubsection{Integral over $J_1\setminus D(s_z)$.}\label{subsubsect:3.3.3}

We split $J_1\setminus D(s_z)$ into four parts:

\bigskip\noindent
(i) $0\le \rho \le \rho_1$ (where $\rho_1$ is a large parameter that will be chosen later);

\bigskip\noindent (ii) $\rho_1 \le \rho \le \delta_3\rho_z$;

\bigskip\noindent (iii) $\delta_3\rho_z \le \rho \le \tfrac32 \rho_z$, $|\rho-\rho_z|\ge \rho_z^{1-\delta_1}$;

\bigskip\noindent  and

\bigskip\noindent (iv) $\rho\ge \tfrac32\rho_z$.

\bigskip\noindent
By Lemma~\ref{lemma:reG-bounded-sector}, the integral over $\rho\in [0, \rho_1]$ is negligible.

\medskip
In the range $\rho_1\le \rho \le \delta_3\rho_z$, we have
\[
\re G(z, \rho e^{{\rm i}\theta_z}) = \rho |\cos\theta_z|\,
\int_\rho^{\rho_z} \frac{\e (u)}{u}\, {\rm d}u + O(\rho)\,.
\]
Consider the function
\[
\ell (\rho) = \rho\, \int_\rho^{\rho_z} \frac{\e (u)}{u}\, {\rm d}u\,.
\]
We have
\[
\ell'(\rho) = \int_\rho^{\rho_z} \frac{\e (u)}{u}\, {\rm d}u - \e (\rho)\,,
\quad \ell''(\rho)
= - \frac{\e (\rho)}{\rho} - \e'(\rho) = -(1+o(1))\frac{\e (\rho)}{\rho}\,.
\]
If $\rho_1$ is chosen sufficiently large, then $\ell''(\rho)$ is negative on
$[\rho_1, \delta_3 \rho_z]$. Hence, $\ell'(\rho)$ decreases. Since
$\delta_3$ is small and fixed,
$\ell'(\delta_3\rho_z) = (1+o(1))(\log\tfrac1{\delta_3}-1)\e(\rho_z) >0$.
Thus, $\ell (\rho)$ attains its maximal value at the end-point $\rho=\delta_3 \rho_z$
where it equals $(1+o(1))\rho_z \e(\rho_z)\, \delta_3\log\tfrac1{\delta_3}$.
This shows that the integral over $[\rho_1, \delta_3\rho_z]$ is negligible, provided that
$\delta_3$ is sufficiently small.

\medskip Now, consider the range $\delta_3\rho_z \le \rho \le \tfrac32 \rho_z$, $|\rho-\rho_z|\ge \rho_z^{1-\delta_1}$. By Lemma~\ref{lemma:ReG-segment},
for $\delta_3\rho_z \le \rho \le \tfrac32 \rho_z$ we have
\[
\frac{\partial^2}{\partial \rho^2}\, \re G(z, \rho e^{{\rm i}\theta_z})
\le -c\,  \frac{\e(\rho_z)}{\rho_z}\,,
\]
whence
\[
\re G(z, \rho e^{{\rm i}\theta_z}) \le \re G(z, s_z) - c\, \frac{\e(\rho_z)}{\rho_z}\,
(\rho-\rho_z)^2\,.
\]
Integrating this over $|\rho-\rho_z|\ge \rho_z^{1-\delta_1}$, we get
a negligible expression.

\medskip The last range to consider is $\tfrac32 \rho_z \le \rho < \infty$.
Here, by Lemma~\ref{lem:ReG-tails},
\[
\frac{\partial}{\partial\rho} \re G(z, \rho e^{{\rm i}\theta_z})
\le - c\e(\rho_z)\,.
\]
Since we already know that
\[
\re G(z, \tfrac32 \rho_z e^{{\rm i}\theta_z}) \le \re G(z, s_z) - c\e(\rho)z\rho_z\,,
\]
we see that the integral over this range is also negligible.

\subsubsection{Integral over $J_2$.}

First, we note that
\[
\re G(z, \rho e^{{\rm i}\theta_z})
- \re G(z, \rho e^{-{\rm i}\theta_z})
= (2+o(1))\rho \e (\rho_z) \theta_z \sin\theta_z
\]
uniformly in $z\in\bar\Omega(\alpha_0-\delta)$, $z\to\infty$.
Thus,
\[
\re G(z, \rho e^{-{\rm i}\theta_z}) \le
\re G(z, \rho e^{{\rm i}\theta_z})
- c \rho \e (\rho_z)\,.
\]
Combining this observation with the estimates of $\re G(z, \rho e^{{\rm i}\theta_z})$
from~\ref{subsubsect:3.3.3}, we readily conclude that the integral over
$J_2$ is negligible as well.

\section{Proof of Theorem~\ref{TheoremE}}

Throughout this section we fix an admissible function $\gamma(s)=L(s)^s$,
which is analytic in the angle $|\arg (s+c_\gamma)|>0$ (with $c_\gamma>0$)
and satisfies conditions
(A), (B), and (C), and assume that
\begin{equation}\label{eq:varepsilon}
\bar\e = \limsup_{\rho\to\infty}\varepsilon(\rho)<2\,.
\end{equation}
We fix positive parameters $\sigma_0$ and $\delta_0$ such that
\[
0<\sigma_0< \min (c_\gamma, 1), \qquad 0 < \delta_0 <
\pi\bigl( \tfrac1{\bar\e} - \tfrac12 \bigr)\,,
\]
and recall that
\[
z E(z) + \frac1{\gamma(0)} = \sum_{n\ge 0} \frac{z^n}{\gamma (n)}\,.
\]
As before, we put $z=re^{{\rm i}\psi}$ with $|\psi|\le \pi$, and
$s= \sigma + {\rm i}t = \rho e^{{\rm i}\theta}$ with $|\theta|<\pi$.

\medskip We will also need the following elementary lower bound for the
function $\gamma$:

\begin{lemma}\label{lemma:gamma}
For $s=\sigma+{\rm i}t$, $\sigma\ge -\sigma_0$, we have
\[
\frac1{|\gamma (s)|} \le C_{\sigma_0}\, \frac{e^{a|t|}}{|\gamma (|\sigma|)|}
\]
with any $a>\tfrac{\pi}2\, \bar\e$. In particular, this holds with some $a<\pi$.
\end{lemma}

\noindent{\em Proof of Lemma~\ref{lemma:gamma}}: We have
\begin{align*}
\log |L(\sigma +{\rm i}t)| &= \int_0^{|\sigma+{\rm i}t|} \frac{\e(u)}u\, {\rm d}u
+ O(1) \\ \\
&\ge \int_0^{|\sigma|} \frac{\e(u)}u\, {\rm d}u + O(1) = \log L(|\sigma|) + O(1)\,,
\end{align*}
and
\[
\bigl| t\, \arg L(\sigma +{\rm i}t) \bigr|
\le \bigl( |\arg (\sigma+{\rm i}t)|\, \e (\sigma +{\rm i} t) |t| + o(|t|)
\le \tfrac{\pi}2\, a |t| + O(1)\,,
\]
with any $a>\bar\e$. This completes the proof of the lemma. \hfill $\Box$

\subsection{Applying the Abel-Plana summation.}

Our starting point is the representation
\begin{align}\label{eq:Abel-Plana}
\nonumber \sum_{n\ge 0} \frac{z^n}{\gamma (n)} &=
\int_{-\sigma_0}^\infty \frac{z^\sigma}{\gamma (s)}\,  {\rm d}\sigma \\
\nonumber \\
&\quad - \frac1{2{\rm i}}\, \int_{-\sigma_0}^{-\sigma_0+{\rm i}\infty}
\frac{z^s}{\gamma (s)} \bigl( \cot (\pi s) + {\rm i} \bigr)\, {\rm d}s
+ \frac1{2{\rm i}}\, \int_{-\sigma_0}^{-\sigma_0-{\rm i}\infty}
\frac{z^s}{\gamma (s)} \bigl( \cot (\pi s) - {\rm i} \bigr)\, {\rm d}s\,.
\end{align}
This is one of the versions of the classical Abel-Plana summation formula.
It holds for any function $F(s)$ holomorphic on $\bigl\{ \re (s) \ge - \sigma_0 \bigr\}$
with convergent series $ \sum_n F(n) $, which satisfies
\[
\lim_{\sigma\to +\infty} |F(\sigma +{\rm i}t)| e^{-b|t|} = 0
\]
with some $b<2\pi$, uniformly in $t\in\mathbb R$. For $F(s) = z^s/\gamma (s)$,
the latter condition immediately follows from Lemma~\ref{lemma:gamma}.

\medskip
Note that the LHS of \eqref{eq:Abel-Plana} is an entire function of $z$,
while the integrals on the RHS are analytic functions in the cut plane
$|\arg (z)|<\pi $ with continuous boundary values on the
upper and lower banks of the cut $\arg (z) = \pm\pi$.

\subsection{Estimating the integrals over vertical lines.}

Here we show that both integrals over vertical lines on the RHS of~\eqref{eq:Abel-Plana}
are $o(1)$ uniformly in $\arg (z)$ when $z\to\infty$, and therefore can be neglected.
Since both estimates follow the same lines, we estimate only the $2$nd integral on the RHS of~\eqref{eq:Abel-Plana}.

\medskip
Noting that
\[
\cot (\pi s) + {\rm i} = 2{\rm i}\, \frac{e^{2\pi {\rm i}s}}{e^{2\pi {\rm i}s}-1}
\]
and recalling that $s=-\sigma_0 +{\rm i}t$ with $0<\sigma_0<1$ and $t\ge 0$, we get
\[
\bigl| \cot (\pi (-\sigma_0 +{\rm i}t) + {\rm i} \bigr|
\le C\, e^{-2\pi t}\,, \qquad t\ge 0\,.
\]
The rest follows from Lemma~\ref{lemma:gamma}:
\[
\Bigl| \int_{-\sigma_0}^{-\sigma_0+{\rm i}\infty}
\frac{z^s}{\gamma (s)} \bigl( \cot (\pi s) + {\rm i} \bigr)\, {\rm d}s \Bigr|
\le C_{\sigma_0}\, \frac{r^{-\sigma_0}}{\gamma (\sigma_0)}\,
\int_0^\infty e^{(-\psi + a - 2\pi) t}\, {\rm d}t\,.
\]
Since $\psi \ge -\pi$ and $a<\pi$, we are done.

\subsection{Estimating the main integral.}

Thus,
\[
\sum_{n\ge 0} \frac{z^n}{\gamma (n)} =
\int_{-\sigma_0}^\infty \frac{z^\sigma}{\gamma (s)}\,  {\rm d}\sigma + o(1)
\]
uniformly in $|\psi|\le \pi$, and
the proof of Theorem~\ref{TheoremE} boils down to estimation of the
integral on the RHS.

Similarly to the proof of Theorem~\ref{TheoremK},
we split the proof into three cases:

\bigskip\noindent (I) $z\in\bar\Omega(\tfrac{\pi}2-\delta_0)$,

\bigskip\noindent (II)
$z\in \bar\Omega(\tfrac{\pi}2+\delta_0)\setminus \Omega(\tfrac{\pi}2-\delta_0)$,

\bigskip\noindent
and

\bigskip\noindent(III)
$z\in \mathbb C\setminus \Omega(\tfrac{\pi}2+\delta_0)$.

\bigskip
In the first two cases, the saddle point $s_z = \rho_z e^{{\rm i}\theta}$
lies in the sectors $|\theta_z|\le \tfrac{\pi}2 - \delta_0$ and
$|\theta_z-\tfrac{\pi}2|\le\delta_0$, correspondingly.
For $z\in\bar\Omega(\tfrac{\pi}2-\delta_0)$, as in the proof of Theorem~\ref{TheoremK}, the
main term comes from integration over a neighbourhood of the saddle point and is given by Lemma~\ref{lemma-s.p.}.
For $z$ close to the boundary of $\Omega(\tfrac{\pi}2)$, the contributions of the saddle point
and of the neighbourhood of the starting point $s=-\sigma_0$ might be of the same order of magnitude.
In the third case, the main term comes from integration over a
neighbourhood of the starting point $s=-\sigma_0$.

\medskip
As in the proof of Theorem~\ref{TheoremK}, we assume
that $\re (z) \ge 0$, that is, $0\le \psi = \arg (z) \le \pi$, and (in cases (I) and (II))
$0\le \theta_z \le \alpha_0-\delta$.

\medskip
As above, we use the notation $G(z, s)=\log \gamma (s) - s \log z$. By $\delta$, $\delta_i$ we denote small positive parameters that remain fixed during our estimates. Most of
these parameters have been already defined in Section~\ref{sect:preliminaries}.

\medskip
In the course of the proof of Theorem~\ref{TheoremE} all expressions that are
\[
o(1) \sqrt{\frac{\rho_z}{\e(\rho_z)}}\, e^{-\re G(z, s_z)} + o(1)
\]
will be called {\em negligible}.

\subsection{Case I: $z\in \bar\Omega(\frac{\pi}{2}-\delta)$.}

In this case, $0\le \theta_z \le \tfrac{\pi}2-\delta_0$, and we will deform the integration
contour to
\[
\Gamma_z = [-\sigma_0, 0] + e^{{\rm i}\theta_z}\,\mathbb R_+\,.
\]
By Lemma~\ref{lemma:gamma}, for $\re (s) \ge 0$ and $\re (z)\ge 0$, we have
\[
\Bigl| \frac{z^s}{\gamma (s)} \Bigr|
\le \frac{C r^\sigma}{\gamma (\sigma)}\,.
\]
For $\sigma\to\infty$ the LHS converges to $0$ faster than exponentially.
This justifies rotation of the integration contour.

\medskip The function $ z^\sigma/\gamma (\sigma) $ remains bounded for
$\sigma\in [-\sigma_0, 0]$ and $|z|\ge 1$. Since the main term of the asymptotics
comes from the integration over
$e^{{\rm i}\theta_z}\mathbb R_+\setminus D(s_z, \rho_z^{1-\delta_1})$ and
grows very fast, we may discard the integration over the segment $ [-\sigma_0, 0] $
and estimate only the integral over
$e^{{\rm i}\theta_z}\mathbb R_+\setminus D(s_z, \rho_z^{1-\delta_1})$.
The estimates we need are practically identical to the ones from~\ref{subsubsect:3.3.3}.
We will not repeat these estimates, only mentioning that
therein we integrated $\exp[\re G(z, \rho e^{{\rm i}\theta_z})]$ along the ray
with strictly negative $\cos(\arg (s))$, while now we
integrate $\exp[-\re G(z, \rho e^{{\rm i}\theta_z})]$ along the ray with
strictly positive $\cos(\arg (s))$.

\subsection{Case II: $z\in \bar\Omega(\frac{\pi}{2}+\delta_0)\setminus \Omega_(\frac{\pi}{2}-\delta_0)$.}

In this case, $\bigl| \theta_z - \tfrac12 \pi \bigr| \le \delta_0$.
We fix an arbitrary positive $\eta < 1/e$, and put
\[
\theta_- = \alpha-\delta, \quad \rho_- = \eta \rho_z, \quad s_- = \rho_- e^{{\rm i}\theta_-},
\quad s_+ = \rho_+ e^{{\rm i}\theta_+} = 2s_z - s_-
\]
(i.e., $s_z$ is the center of the interval $[s_-, s_+]$), and choose $t_0$ so that
$\arg (-\sigma_0 + {\rm i} t_0) = \theta_-$. 
\begin{figure}[h]
 \includegraphics[scale=0.8]{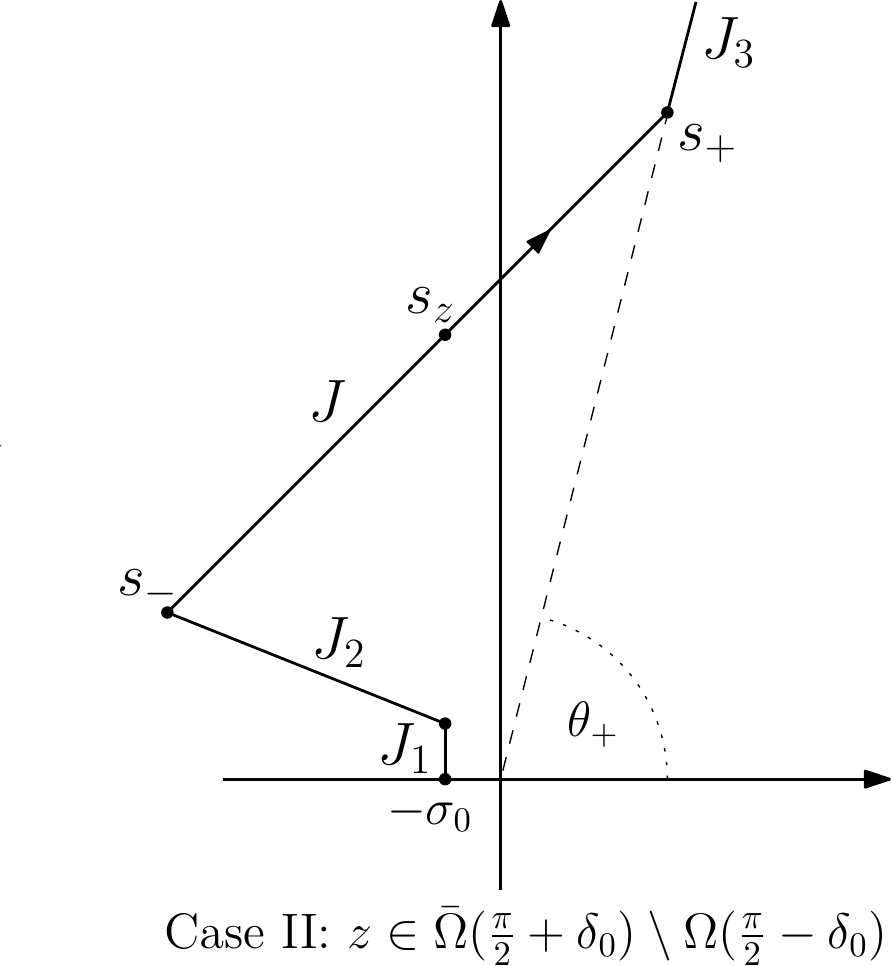}
   \caption{$\Gamma_{z} $}
\end{figure}
Then, we deform the integration contour
to the union of three segments and a ray
\begin{align*}
\Gamma_z &= J_1 + J_2 + J + J_3 \\ \\
&=
[-\sigma_0, -\sigma_0 +{\rm i} t_0] +
[-\sigma_0 +{\rm i} t_0, s_-] +
[s_-, s_+] + e^{{\rm i}\theta_+} [\rho_+, \infty)\,.
\end{align*}
We note that if the parameter $\delta_0$ is sufficiently small, then
$0<\theta_+<\min\bigl( \tfrac{\pi}2, \theta_z\bigr)$. This justifies
the deformation of the contour.

\medskip
Next, we note that the segment $J$ traverses the saddle-point $s_z$ in the direction
$ \arg (s-s_z)= \tfrac{\pi}2-c\eta$ provided that $\delta_0$ is sufficiently small.
To apply Lemma~\ref{lemma-s.p.}, we need to traverse the saddle point in the direction
$ |\arg (s-s_z)-\tfrac12 \theta_z| \le \tfrac{\pi}4-\delta_2$; taking into account
that $\bigl| \theta_z - \tfrac12 \pi \bigr| \le \delta_0$, this means that
$ \tfrac12 \delta_0 + \delta_2 \le \arg (s-s_z) \le \tfrac{\pi}2 -
(\tfrac12 \delta_0 + \delta_2)$. If $\delta_0$ and $\delta_2$ are
sufficiently small, the direction $ \tfrac{\pi}2-c\eta $ lies within
this range, and therefore, Lemma~\ref{lemma-s.p.} is applicable. It tells
us that the integral over $J\cap D(s_z, \rho_z^{1-\delta_1})$ equals
\[
(1+o(1))\, \sqrt{\frac{2\pi s_z}{\e (s_z)}}\, \exp \bigl[ s_z\e (s_z) \bigr]\,.
\]
It remains to see that the contributions of the integrals over
$J\setminus D(s_z, \rho_z^{1-\delta_1})$, $J_1$, $J_2$, and $J_3$ are
all negligible.

\subsubsection{Integral over $J_1$.}

For $s\in J_1$, we have $|z^s/\gamma (s)|\le C r^{-\sigma_0}$. Hence,
the integral over $J_1$ is bounded by $ C r^{-\sigma_0} $ and 
can be neglected.

\subsubsection{Integral over $J_2$.}\label{subsubsect:EII-J2}

As in the previous case, for any given $\rho_1>\rho_0$ (independent of $z$) and
for $s\in [\rho_0 e^{{\rm i}\theta_-}, \rho_1 e^{{\rm i}\theta_-}]$, we have
$|z^s/\gamma (s)|\le C r^{-\sigma_0}$. Hence, integrating over $J_2$, we can
integrate only over $\rho_1 \le \rho \le \eta \rho_z$.
Then, by~\eqref{eq:ReG-asymptotics},
\[
- \re G(z, \rho e^{{\rm i}\theta_-}) = - \rho |\cos\theta_-|\,
\int_{\rho}^{\rho_z} \frac{\e (u)}{u}\, {\rm d}u + O(\rho)
\]
uniformly in $z\in\bar\Omega (\alpha-\delta)$.

We claim that, {\em given positive $A$ and $\eta<1/e$, we have
\[
\rho\, \int_\rho^{\rho_z} \frac{\e (u)}{u}\, {\rm d}u > A \log \rho\,,
\qquad \rho_1 \le \rho \le \eta\rho_z\,,
\]
provided that $\rho_1$ is sufficiently large}. Indeed, consider the
function
\[
\ell (\rho) = \rho\, \int_\rho^{\rho_z} \frac{\e (u)}{u}\, {\rm d}u - A \log \rho\,.
\]
We have
\[
\ell'(\rho) = \int_\rho^{\rho_z} \frac{\e (u)}{u}\, {\rm d}u
- \e(\rho) - \frac{A}{\rho}\,,
\]
and
\[
\ell''(\rho) = - \frac{\e(\rho)}{\rho} - \e'(\rho) + \frac{A}{\rho^2}
= - (1+o(1))\, \frac{\e (\rho)}{\rho} < 0\,,
\]
whenever $\rho_1$ is sufficiently large. Therefore, the function $\ell'(\rho)$
decays on $[\rho_1, +\infty)$. Noting that
\[
\ell'(\eta \rho_z) = (1+o(1)) \e (\rho_z) \bigl( \log\tfrac{1}{\eta} - 1\bigr)
\]
and recalling that $\eta<1/e$, we conclude that $\ell'>0$ on $[\rho_1, \eta \rho_z]$,
i.e., $\ell $ increases therein. Therefore,
\[
\ell (\rho) \ge \ell (\rho_1)
= \rho_1\, \int_{\rho_1}^{\rho_z} \frac{\e (u)}{u}\, {\rm d}u - \frac{A}{\rho_1} > 0\,,
\]
provided that $z$ is large enough. This proves the claim.

\medskip This claim immediately yields that
\[
- \re G(z, \rho e^{{\rm i}\theta_-}) \le -2 \log\rho\,, \qquad
\rho_1 \le \rho \le \eta\rho_z\,.
\]
Hence, the integral we are estimating  does not exceed
\[
\int_{\rho_1}^\infty\, \frac{{\rm d}\rho}{\rho^2} = \frac1{\rho_1}\,.
\]
Since we can choose $\rho_1$ as large as we need, the integral over $J_2$ is
negligible.

\subsubsection{Integral over $J\setminus D(s_z, \rho_z^{1-\delta_1})$.}

For $s\in J$, we have $\arg (s-s_z)=\tfrac{\pi}2-c\eta$ and
$\arg (s)\in (\theta_+, \theta_-)$. Therefore,
\[
\pi - 2c\eta - \theta_+ \le 2\arg (s-s_z) - \arg (s) \le
\pi - 2 c \eta - \theta_+\,.
\]
So, if $\eta$ is chosen sufficiently small, $ \cos( 2\arg (s-s_z) - \arg (s)) \ge c >0$
for every $s\in J$. Then, by Lemma~\ref{lemma:ReG-segment}, the function
$\re G(z,\, \cdot \,)$ is concave on $J$, and moreover
\[
\frac{\partial^2}{\partial t^2} \bigl[ - \re G(z, s_z + t e^{{\rm i}\phi}) \bigr]
\le - c \rho_z \e (\rho_z)
\]
whenever $s_z+te^{{\rm i}\phi}\in J$. The rest of the argument is the same as
in~\ref{subsubsect:KII-segment} and we skip it.

\subsubsection{Integral over $J_3$.}

We skip the estimate since it follows the same lines as the one in~\ref{subsubsect:J_1-J_3}.

\subsection{Case III: $z\in \mathbb C\setminus \bar\Omega (\frac{\pi}{2}+\delta_0)$.}

In this case, the saddle-point $s_z$ may not exist (more precisely, it may live on another sheet
of the Riemann surface of $\log z$). Nevertheless, given $z=r e^{i\psi}$, we define a positive value $\rho_z$ by
equation
\[
\int_0^{\rho_z} \frac{\varepsilon(u)}{u}\, {\rm d}u = \log r\,.
\]
\begin{figure}[ht]
\includegraphics[scale=0.75]{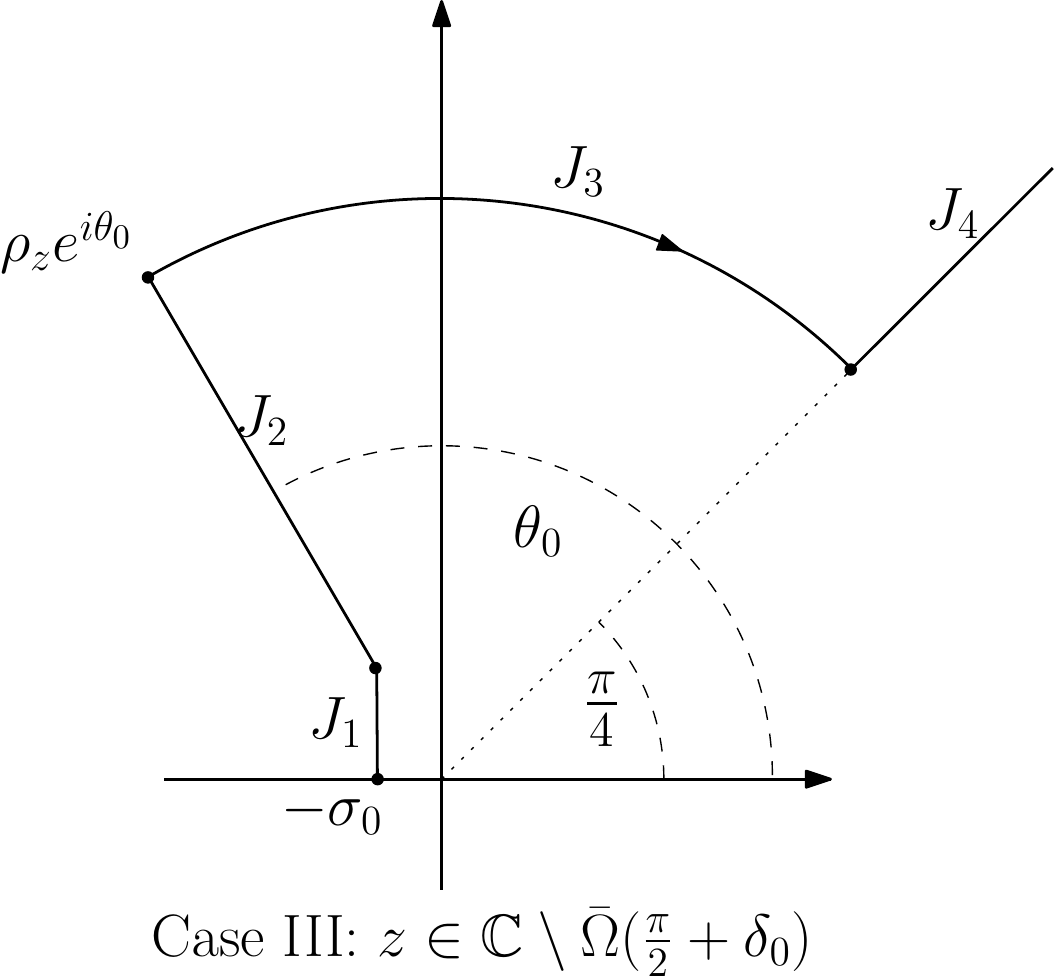}
\caption{$\Gamma_{z} $}
\end{figure}
Then, we choose $t_0>0$ so that $\arg (-\sigma_0 +{\rm i} t_0)$,
put $\theta_0=\tfrac{\pi + \delta_0}2$, $\rho_0=|\sigma_0+{\rm i} t_0|$,
and deform the contour to
\begin{align*}
\Gamma_z &= J_1 + J_2 - J_3 + J_4 \\ \\
&= [-\sigma_0, - \sigma_0+{\rm i} t_0] +
[\rho_0 e^{{\rm i} \theta_0}, \rho_z e^{{\rm i} \theta_0}]
- \bigl\{ \rho_z e^{{\rm i}\theta}\colon \tfrac{\pi}4 \le\theta\le \theta_0 \bigr\}
+ e^{{\rm i}\pi/4}[\rho_z, +\infty)\,.
\end{align*}
As in the previous cases, the integral over $J_1$ is bounded by
$C r^{-\sigma_0}$, and we need to estimate the other three integrals.

\subsubsection{Integral over $J_2$.}

For sufficiently large $z\in \mathbb C\setminus \bar\Omega (\frac{\pi}{2}+\delta_0 )$,
we have
\[
\psi = \arg (z) \ge \bigl( \tfrac{\pi}2+\tfrac34 \delta_0 \bigr) \e(\rho_z),
\]
while for $s\in J_2$, $ \theta_0 = \arg (s) = \tfrac{\pi}2+\tfrac12 \delta_0 $.
Therefore,
\begin{align}\label{eq:-reG}
- \re G(z, \rho e^{{\rm i}\theta_0}) &\le - \rho|\cos \theta_0|\,
\int_\rho^{\rho_z} \frac{\e(u)}u\, {\rm d} u
- \rho\bigl( \psi - \theta_0 \e(\rho) \bigr)\sin\theta_0  + \rho o(\e (\rho))
\\ \nonumber \\ \nonumber
&\le  - \rho\,
\int_\rho^{\rho_z} \Bigl( \frac{\e(u)}u |\cos\theta_0|
+ \bigl( \tfrac{\pi}2+\tfrac12 \delta_0 \bigr) \e'(u) \sin\theta_0 \Bigr)\,
{\rm d} u - \tfrac14 \delta_0  \rho \e (\rho_z)\sin\theta_0  + \rho o(\e (\rho))
\\ \nonumber \\ \nonumber
&\le - c\delta_0 \rho \Bigl[ \int_\rho^{\rho_z} \frac{\e(u)}u\, {\rm d} u  + \e (\rho_z) \Bigr]\,.
\end{align}
Then, as in the estimate of the similar integral in~\ref{subsubsect:EII-J2},
we split the integral into three pieces: (a) $\rho_0 \le \rho \le \rho_1$ with sufficiently
large $\rho_1$, (b) $\rho_1 \le \rho \le \eta \rho_z$, and (c) $\eta\rho_z \le \rho \le \rho_z$.

\medskip\noindent (a)
For $\rho_0 \le \rho \le \rho_1$, as before,
we have $ | z^s/\gamma (s) | \le C r^{-\sigma_0} $, which does the job.

\medskip\noindent (b)
Then, similarly to~\ref{subsubsect:EII-J2}, we choose sufficiently large $\rho_1$
and $\eta$ much smaller than $\delta_0$, so that, for $\rho_1 \le \rho \le \eta \rho_z$,
\[
c\delta_0 \rho \int_\rho^{\rho_z} \frac{\e(u)}u\, {\rm d} u > 2 \log r\,.
\]
This bounds the integral over the second piece by $1/\rho_1$.

\medskip\noindent (c)
At last, on the third piece, we use that
\[
- \re G(z, \rho e^{{\rm i}\theta_0}) \le - c\delta_0\rho\e (\rho_z)
\le -c \delta_0 \eta \rho_z \e(\rho_z)\,,
\]
which is, by far, more than we need.

\subsubsection{Integral over $J_3$.}

Using estimate~\eqref{eq:-reG}, and recalling that
\[
\psi \ge \bigl( \tfrac{\pi}2+\tfrac34 \delta_0 \bigr) \e(\rho_z),
\]
while
$0 \le \theta\le \tfrac{\pi}2+\tfrac34 \delta_0 $, we get
\[
-\re G(z, \rho_z e^{{\rm i}\theta}) \le -\rho_z (\psi - \theta\e (\rho_z))
+ \rho_z o(\e (\rho_z)) \le -c\delta_0 \rho_z \e(\rho_z)\,,
\]
which immediately yields that the integral over $J_3$ decays fast to zero.

\subsubsection{Integral over $J_4$.}

Here using once again estimate~\eqref{eq:-reG} and recalling that
$\psi < -\tfrac{\pi}2\, \e(\rho_z)$, $\theta=\tfrac{\pi}4$, and
$\rho\ge\rho_z$, we get
\begin{align*}
-\re G(z, \rho e^{{\rm i}\pi/4}) &\le
- \frac{\rho}{\sqrt{2}}\, \int_{\rho_z}^\rho
\Bigl( \frac{\e (u)}u - \frac{\pi}4\, \e'(u) \Bigr)\, {\rm d}u
-\frac{\pi}4\, \frac{\rho}{\sqrt 2}\, \e(\rho_z) + \rho \e (\rho) \\ \\
&\le -c\rho \Bigl[  \int_{\rho_z}^\rho \frac{\e (u)}u\, {\rm d}u
+ \e (\rho_z) \Bigr] \\ \\
&\le -c\rho \e (\rho_z)\,,
\end{align*}
and then,
\[
\int_{\rho_z}^\infty e^{-c\rho \e (\rho_z)}\, {\rm d}\rho
\le \frac{C}{\e (\rho_z)}\, e^{-c\rho_z \e (\rho_z)}\,,
\]
which decays very fast to zero, uniformly in
$z\in \mathbb C\setminus \bar\Omega (\frac{\pi}{2}+\delta_0 )$, $z\to\infty$.
This completes the proof of Theorem~\ref{TheoremE}. \hfill $\Box$

\begin{appendices}

\section{Proof of Theorem~3}

\subsection{Slowly varying function.}
We will need several well-known proprieties of slowly-varying functions, which we summarize in the following lemma:

\begin{lemma}\label{lemma-A1}
Suppose that the function $\varepsilon\in C^1[0, \infty)$
is such that $\rho \varepsilon'(\rho) = o(\varepsilon(\rho))$, $\rho\to\infty$.
Then,

\medskip\noindent{\rm 1.\ }
For any interval $[a, b]\subset (0, \infty)$,
$\displaystyle \lim_{\rho\to\infty}\, \sup_{\lambda\in[a,b]}\, \frac{\varepsilon(\lambda \rho)}{\varepsilon(\rho)}=1$.

\medskip\noindent{\rm 2.\ }
For any $\delta>0$, the function $\rho\mapsto \rho^\delta \varepsilon(\rho)$ is eventually increasing and the function $\rho\mapsto \rho^{-\delta} \varepsilon(\rho)$ is eventually decreasing.

\bigskip\noindent{\rm 3.\ }
If  $I$ is an interval,  $m\in C(\mathbb{R}_+ \times I)$ and $\delta>0$,  such that $\int_0 t^{-\delta}|m(t,x)|\,{\rm d}t <\infty $ and $\int^\infty t^{\delta}|m(t,x)|\, {\rm d}t <\infty,$ for all $x\in I$, then
\[\lim_{\rho\to\infty} \int_0^\infty\frac{\varepsilon(\rho t)}{\varepsilon(\rho)}\, m(t,x)\, {\rm d}t= \int_0^\infty m(t,x)\, {\rm d}t\]
locally uniformly in $I$.
\end{lemma}

\medskip
The proofs of assertions 1--3 can be found in~\cite{Bingham}, Lemmas/Theorems 1.3.1, 1.5.5 and 4.5.2 correspondingly.

\subsection{Two lemmas that yield Theorem~3.}
We fix the function $\ell$ that satisfies  assumptions of Theorem~3, that is, $\ell:[0,\infty)\to(0,\infty)$ is an unboundedly increasing $C^1$--function such that the function
 	\[ \rho\mapsto \rho\frac{\ell'(\rho)}{\ell(\rho)} \]
 	is slowly varying and bounded on $[0,\infty)$,
 	and put
\begin{equation}
\gamma(s)=\exp\left(s^2\int_0^\infty \frac{\ell'(u)}{\ell(u)}\cdot \frac{{\rm d}u}{s+u} \right),\quad |\arg (s)|<\pi,\label{newGamma}\,.
\end{equation}
Then,
 \[ \log L(s)=\frac1{s}\, \log\gamma(s)=s\int_0^\infty \frac{\ell'(u)}{\ell(u)} \cdot\frac{{\rm d}u}{s+u},\quad |\arg (s)|<\pi, \]
 and
\[
\varepsilon (s) = s\,\frac{L'(s)}{L(s)}=s\int_0^\infty u\,\frac{\ell'(u)}{\ell(u)}\, \frac{{\rm d}u}{(s+u)^2}, \quad |\arg(s)|<\pi\,.
\]	
The following two lemmas immediately yield Theorem~3.

\begin{lemma}\label{lemma-A2}
The function $\gamma$ defined by \eqref{newGamma} is admissible.
\end{lemma}

\begin{lemma}\label{lemma-A3} We have

\medskip\noindent{\rm 1.\ }
$\displaystyle \lim_{\rho\to\infty}\frac{\log L(\rho)}{\log \ell(\rho)}=1. $

\medskip\noindent{\rm 2.\ }
If in addition there exists the limit $\displaystyle\lim_{\rho\to\infty} \rho\frac{\ell'(\rho)}{\ell(\rho)}$, then,  $\displaystyle \lim_{\rho\to\infty}\frac{ \ell(\rho)}{ L(\rho)}=\ell(0). $
\end{lemma}

\subsection{Proof of Lemma~\ref{lemma-A2}.}
By our assumption,  $\tfrac{\ell'(u)}{\ell(u)}=O\left(\frac{1}{u}\right)$
as $u\to\infty$. Therefore, the integral in the definition of the function
$\gamma$ is absolutely and locally uniformly convergent in $\{s\colon |\arg(s)|<\pi\}$,
and therefore, the function $\gamma$ is analytic and non-vanishing therein. It is
easy to see that positivity of $\ell'$ yields continuity of $\tfrac1{\gamma}(s)$ at $s=0$.

By Lemma~\ref{lemma-A3}, assertion~3, applied with
\[
m(t, \theta) = \frac{e^{{\rm i}\theta}}{(e^{{\rm i}\theta}+t)^2},
\quad e^{{\rm i}\theta}\, \int_0^\infty \frac{{\rm d}t}{(e^{{\rm i}\theta}+t)^2} = 1\,,
\]
we get
 	\begin{equation}
 	\varepsilon (s) =\left(1+o(1)\right) \rho\,\frac{\ell'(\rho)}{\ell(\rho)},\quad s=\rho e^{{\rm i}\theta}, \rho\to\infty, \label{ToproveC}
 	\end{equation}
 	uniformly in any angle $|\arg(s)|\le\pi-\delta$. This gives us the properties (A) and (C) in the definition of admissible functions.
 	
In order to show that
the property (B) also holds, we differentiate under the integral sign once again, and obtain
 \[  \rho\varepsilon'(\rho)=-\rho\int_0^\infty u\frac{\ell'(u)}{\ell(u)} \cdot\frac{\rho-u}{(u+\rho)^3}\,{\rm d}u\,.\]
 Since,
 \[
 \int_0^\infty \frac{\rho-u}{(\rho+u)^3}{\rm d}u=\int_0^\infty \frac{1-t}{(1+t)^3}\, {\rm d}t = 0\,,  \]
applying again Lemma~\ref{lemma-A3}, this time with $m(t) = (1-t)/(1+t)^3$, we find that
 \[
 \rho\varepsilon'(\rho)=o(1)\cdot \rho\frac{\ell'(\rho)}{\ell(\rho)},\quad \rho\to\infty\,. \]
By~\eqref{ToproveC}, this gives us the property (B). \qed

\subsection{Proof of Lemma~\ref{lemma-A3}.}
  	
\subsubsection{Proof of Part 1.}
 	Integration by part yields
 \begin{multline}
 \log L(\rho)
 = \rho\int_0^\infty \frac{\ell'(u)}{\ell(u)} \cdot\frac{{\rm d}u}{\rho+u}
 = \rho\int_0^\infty  \log \frac{\ell(u)}{\ell(0)}\cdot\frac{{\rm d}u}{(\rho+u)^2} \\
 = \int_0^\infty  \log \frac{\ell(\rho t)}{\ell(0)}\cdot\frac{{\rm d}t}{(1+t)^2}\,.
 \label{ellafterint}
 \end{multline}	
 	 	Since
 	 	\[
 \lim_{u\to\infty }u\frac{\left(\log \ell(u)\right)'}{\log \ell(u)} =
 \lim_{u\to\infty }u\frac{\ell'(u)}{\ell(u)}\cdot \frac{1}{\log \ell(u)}=0\,,
 \]
 	 	the function $u\mapsto \log(\ell(u)/\ell(0))$ is slowly varying. Therefore, by Lemma~A.1, 
 assertion~3, applied with $m(t)=(1+t)^{-2}$, we get
 	 	\[ \lim_{\rho\to\infty}\frac{\log L(\rho)}{\log \ell(\rho)-\log \ell(0)}=\int_0^\infty \frac{{\rm d}t}{(1+t)^2}=1.  \]
 	 	Taking into account that $\displaystyle \lim_{\rho\to\infty}\ell(\rho)=\infty$, we conclude 
 assertion~1.

\subsubsection{Part 2.}
By \eqref{ellafterint},	
  	\begin{align*}
  	\log L(\rho)-\log \frac{\ell(\rho)}{\ell(0)}&=\int_0^\infty  \log \frac{\ell(\rho t)}{\ell(0)}\cdot\frac{{\rm d}t}{(1+t)^2}- \log\frac{\ell(\rho)}{\ell(0)}\int_0^{\infty}  \frac{{\rm d}u} {(1+t)^2} \\&=
  	\int_0^\infty \log \frac{\ell(\rho t)}{\ell(\rho)} \cdot \frac{1}{(1+t)^2}\, {\rm d}t.
  	\end{align*}
  We need to show that the integral on the RHS tends to $0$ as $\rho\to\infty$.

  	We fix $\lambda>1$, split the integral into three parts
  	\[
  \int_0^\infty \log \frac{\ell(\rho t)}{\ell(\rho)} \cdot \frac{1}{(1+t)^2}\, {\rm d}u = \int_0^{\lambda^{-1}} +\int_{\lambda^{-1}}^{\lambda}+\int_{\lambda}^{\infty} = I+II+III\,,
  	\]
  	and estimate each of the three integrals separately.
  	
  By our assumption,  $\tfrac{\ell'(u)}{\ell(u)}=O\left(\frac{1}{u}\right)$
  	 as $u\to\infty$. Hence,
  	 \[ \left| \log \frac{\ell(\rho t)}{\ell(\rho)}\right|=\left|\int_\rho^{\rho t} \frac{\ell'(u)}{\ell(u)}\, {\rm d}u \right|\leq C|\log t|\,. \]
Therefore,
  	\[
  	|I| \leq C\int_0^{\lambda^{-1}} \log \frac{\ell (\rho)}{\ell(\rho t )}\,  {\rm d}t \leq  C \int_0^{\lambda^{-1}} |\log t|\, {\rm d}t\leq C\frac{\log \lambda}{\lambda}\,.
  	\]
Similarly,
  	\[
 |III|  \leq C \int_{\lambda}^{\infty}\frac{\log t}{(t+1)^2}\, {\rm d}t
 \leq C \frac{\log \lambda} {\lambda}\,.
  	\]
  	
Letting $a=\lim_{\rho\to\infty }\rho\frac{\ell'(\rho)}{\ell(\rho)}$, we see that
\[ \lim_{\rho\to\infty} \log \frac{\ell (t\rho )}{\ell(\rho )}=\lim_{\rho\to\infty}\int_{\rho}^{\rho t}\frac {\ell'(u)}{\ell(u)}\, {\rm d}u=a\log t\,. \]
uniformly in $t\in [\lambda^{-1},\lambda]$. Therefore,
\[\lim_{\rho\to\infty } II= \lim_{\rho\to\infty } a\int_{\lambda^{-1}}^{\lambda}\frac{\log t}{(t+1)^2}\, {\rm d}t=0\,. \]
  	Letting $\lambda \to \infty$, we obtain
  	\[
  	\lim_{\rho \to \infty}\left[  \log L(\rho)-\log \frac{\ell(\rho)}{\ell(0)}\right]=0.
  	\]
  	This finishes the proof of Lemma~\ref{lemma-A3}, and hence, of Theorem~\ref{thm3}.\qed

\section{Admissible functions of positive type}

Following Beurling, we say that a function $\gamma$ is \textit{of positive type} if it can be represented in the form
\[ \gamma(s)=\exp\left(A+Bs+(s-a)^2\int_0^\infty \frac{{\rm d}\mu(s)}{u+s}\right)\,, \]
with real constants $A,B$ and $a$, and with a non-decreasing function $\mu$ such that
\[\int_0^\infty \frac{d\mu(s)}{u+1}<\infty.  \]

\subsection{} Note that the functions $\Gamma(s)$ and $\log_k^{\beta s}(s+c)$, $\beta>0$, $c>0$ and sufficiently big, are of positive type (here, as before, $\log_k$ is the $k$th
iterate of the logarithmic function). In the first case, this follows from the classical representation
\[ \log \Gamma(s)=\log s-as+s^2\int_0^\infty \frac{[u]}{u^2}\frac{{\rm d}u}{u+s}\,, \]
where $[u]$ denotes the integer part of $u$.
To see this in the second case, we put
\[ F(s)=\frac{\log_k^b(s+c)-\log_k^b(c)}{s},\quad |\arg z|<\pi \]
and apply the Cauchy formula
\[
F(s)=\frac{1}{2\pi i }\int_{\Gamma_{\delta,R}}\frac{F(w)dw}{w-s},\quad 0<\delta\ll 1,
\ R\gg1,
\]
where the contour $\Gamma_{\delta,R}$ is defined in Figure 2.
\begin{figure}[h]
	\centering
	\includegraphics[scale=0.6]{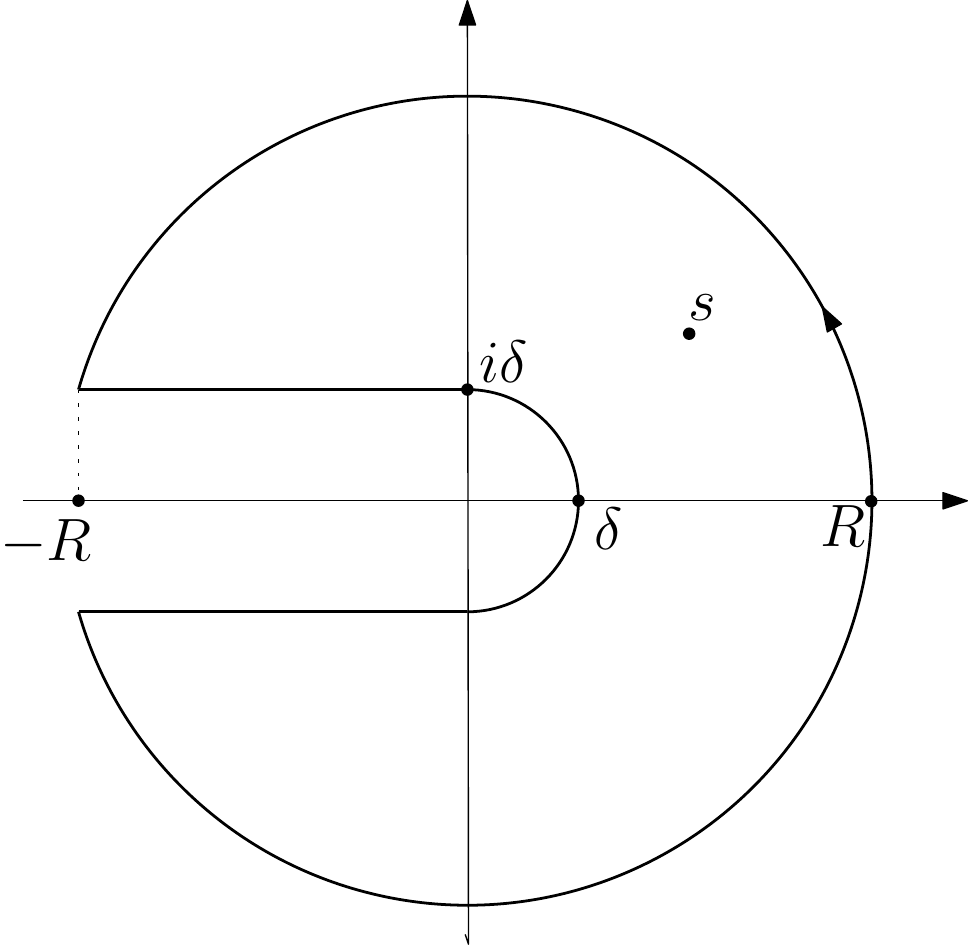}
	\caption{Contour  $\Gamma_{\delta,R}$}
\end{figure}
Then, letting $\delta\to 0$ and $R\to\infty$, we obtain the representation
\[ \log_k^b(s+c)=\log_k^b(c)+s\int_0^\infty \frac{\lambda(u)}{u}\frac{{\rm d}u}{u+s}, \]
	with
	\[ \lambda(u)=\frac{1}{2\pi {\rm i}} \lim_{\epsilon\downarrow 0}\left[\log_k^b(-u+c+{\rm i}\epsilon)-\log_k^b(-u+c-{\rm i}\epsilon)\right]. \]
If the constant $c$ is big enough, the function $\log_k^b(u+c)$ is real on the positive half-line. Hence, 
by the Schwarz reflection principle, its jump on the negative ray is purely imaginary. Then, it is not difficult to see that $\lambda(u)$ is positive, provided that the constant $c$ is big.

\subsection{}

It is easy to check that  functions constructed from these two examples using the rules described in 1.5.1, 1.5.4 and 1.5.5 are also functions of positive type.

\subsection{} We finish this discussion with the statement of Beurling's theorem~\cite{Beurling2, Beurling}:

\medskip\noindent{\bf Theorem (Beurling).}
{\em	Every function $\gamma$ of positive type is Mellin-positive definite, that is, represented by an absolutely convergent Stieltjes integral
	\begin{equation}
	\gamma(s)=\int_0^\infty t^{s-1}\, {\rm d}\nu(t)\label{MPD}
	\end{equation}
	with a non--decreasing function $\nu$ on $[0,\infty)$.
}

\medskip
It worth mentioning that this theorem can be also deduced from results presented in Berg's work~\cite{Berg1,Berg2}.

\medskip
Note that, up to a normalization, the function
$\nu$ in \eqref{MPD} is unique and can be recovered from $\gamma$ by the inverse Mellin transform \cite[\S6.9 ]{Widder}:
\[
\nu(t) = \lim_{T\to+\infty}\frac{{\rm i}}{2\pi}
\int_{c-{\rm i}T}^{c+{\rm i}T}\gamma(s)t^{-z}\frac{{\rm d}s}{s}\,,
\quad c>0\,.
\]
Juxtaposing this with the representation~\eqref{eq:K_Mellin} of the function $K$, we conclude that $K=\nu'$.

\subsection{} Finally, it worth mentioning that, the function $K$ is the unique solution to the Stieltjes moment problem with the moments $(\gamma(n+1))_{n\geq 0}$ whenever
\begin{equation}
\limsup_{\rho\to\infty}\varepsilon(\rho)<2.\label{varepsilon2}
\end{equation}
where, as before, $\varepsilon(\rho)=\rho\frac{L'(\rho)}{L(\rho)}$, in particular,
whenever $\varepsilon(\rho)\to 0$ as $\rho\to \infty$. This immediately follows from the Carleman sufficient condition of determinacy: condition \eqref{varepsilon2} readily yields that $L(\rho)=O(\rho^c),$ $\rho\to\infty$ with some $c<2$, which, in turn, yields divergence of the series
\[ \sum_{n\geq 0}\gamma(n+1)^{-1/(2n)}\geq \sum_{n\geq 1}\frac{1}{\sqrt{L(n)}}=+\infty. \]

\end{appendices}

\end{document}